\documentclass[12pt]{amsart}
\usepackage[nobysame,abbrev,alphabetic]{amsrefs}
\usepackage{amssymb}
\usepackage[all]{xy}

\DeclareMathOperator{\Aut}{Aut}

\DeclareMathOperator{\Gal}{Gal}
\DeclareMathOperator{\Hom}{Hom}
\DeclareMathOperator{\id}{id}
\DeclareMathOperator{\Img}{Im}
\DeclareMathOperator{\Inf}{Inf}

\DeclareMathOperator{\Ker}{Ker}

\DeclareMathOperator{\Res}{Res}
\DeclareMathOperator{\rk}{rk}
\DeclareMathOperator{\syml}{sl}
\DeclareMathOperator{\GL}{GL}

\begin{document}

\newtheorem{thm}{Theorem}[section]
\newtheorem{cor}[thm]{Corollary}
\newtheorem{lem}[thm]{Lemma}
\newtheorem{prop}[thm]{Proposition}
\newtheorem{exam}[thm]{Example}
\newtheorem{examples}[thm]{Examples}
\newtheorem{rem}[thm]{Remark}
\newtheorem{case}{\sl Case}
\newtheorem{conj}[thm]{Conjecture}
\newtheorem*{thmA}{Theorem A}
\newtheorem*{thmB}{Theorem B}
\swapnumbers
\newtheorem{rems}[thm]{Remarks}
\newtheorem*{acknowledgment}{Acknowledgment}
\numberwithin{equation}{section}

\newcommand{\dirlim}{\varinjlim}
\newcommand{\discup}{\ \ensuremath{\mathaccent\cdot\cup}}
\newcommand{\gr}{\mathrm{gr}}
\newcommand{\nek}{,\ldots,}
\newcommand{\inv}{^{-1}}
\newcommand{\isom}{\cong}
\newcommand{\ndiv}{\hbox{$\,\not|\,$}}
\newcommand{\proj}{\mathrm{proj}}
\newcommand{\sep}{\mathrm{sep}}
\newcommand{\tensor}{\otimes}

\newcommand{\alp}{\alpha}
\newcommand{\gam}{\gamma}
\newcommand{\Gam}{\Gamma}
\newcommand{\del}{\delta}
\newcommand{\Del}{\Delta}
\newcommand{\eps}{\epsilon}
\newcommand{\lam}{\lambda}
\newcommand{\Lam}{\Lambda}
\newcommand{\sig}{\sigma}
\newcommand{\Sig}{\Sigma}
\newcommand{\bfA}{\mathbf{A}}
\newcommand{\bfB}{\mathbf{B}}
\newcommand{\bfC}{\mathbf{C}}
\newcommand{\bfF}{\mathbf{F}}
\newcommand{\bfP}{\mathbf{P}}
\newcommand{\bfQ}{\mathbf{Q}}
\newcommand{\bfR}{\mathbf{R}}
\newcommand{\bfS}{\mathbf{S}}
\newcommand{\bfT}{\mathbf{T}}
\newcommand{\bfZ}{\mathbf{Z}}
\newcommand{\dbA}{\mathbb{A}}
\newcommand{\dbC}{\mathbb{C}}
\newcommand{\dbF}{\mathbb{F}}
\newcommand{\dbG}{\mathbb{G}}
\newcommand{\dbI}{\mathbb{I}}
\newcommand{\dbK}{\mathbb{K}}
\newcommand{\dbN}{\mathbb{N}}
\newcommand{\dbP}{\mathbb{P}}
\newcommand{\dbQ}{\mathbb{Q}}
\newcommand{\dbR}{\mathbb{R}}
\newcommand{\dbU}{\mathbb{U}}
\newcommand{\dbV}{\mathbb{V}}
\newcommand{\dbZ}{\mathbb{Z}}
\newcommand{\grf}{\mathfrak{f}}
\newcommand{\gra}{\mathfrak{a}}
\newcommand{\grA}{\mathfrak{A}}
\newcommand{\grB}{\mathfrak{B}}
\newcommand{\grh}{\mathfrak{h}}
\newcommand{\grH}{\mathfrak{H}}
\newcommand{\grI}{\mathfrak{I}}
\newcommand{\grL}{\mathfrak{L}}
\newcommand{\grm}{\mathfrak{m}}
\newcommand{\grp}{\mathfrak{p}}
\newcommand{\grq}{\mathfrak{q}}
\newcommand{\grr}{\mathfrak{r}}
\newcommand{\grR}{\mathfrak{R}}
\newcommand{\grU}{\mathfrak{U}}
\newcommand{\grZ}{\mathfrak{Z}}
\newcommand{\calA}{\mathcal{A}}
\newcommand{\calB}{\mathcal{B}}
\newcommand{\calC}{\mathcal{C}}
\newcommand{\calE}{\mathcal{E}}
\newcommand{\calG}{\mathcal{G}}
\newcommand{\calH}{\mathcal{H}}
\newcommand{\calJ}{\mathcal{J}}
\newcommand{\calK}{\mathcal{K}}
\newcommand{\calL}{\mathcal{L}}
\newcommand{\calR}{\mathcal{R}}
\newcommand{\calW}{\mathcal{W}}
\newcommand{\calU}{\mathcal{U}}
\newcommand{\calV}{\mathcal{V}}
\newcommand{\calZ}{\mathcal{Z}}

\newcommand{\cc}{\mathfrak c}
\newcommand{\dd}{\mathfrak d}

\title[Symbol length]{The symbol length for elementary type pro-$P$ groups and Massey products}

\author{Ido Efrat}

\address{Earl Katz Family Chair in Pure Mathematics,
Department of Mathematics, Ben-Gurion University of the Negev,
P.O. Box 653, Be'er-Sheva 8410501, Israel}

\email{efrat@bgu.ac.il}

\begin{abstract}
For a prime number $p$ and an integer $m\geq2$, we prove that the symbol length of all elements of $m$-fold Massey products in $H^2(G,\dbF_p)$, for pro-$p$ groups $G$ of elementary type, is bounded by $(m^2/4)+m$.
Assuming the Elementary Type Conjecture, this applies to all finitely generated maximal pro-$p$ Galois groups $G=G_F(p)$ of fields $F$ which contain a root of unity of order $p$.
More generally, we provide such a uniform bound for the symbol length of all pullbacks $\rho^*(\bar\omega)$ of a given cohomology element $\bar\omega\in H^n(\bar G,\dbF_p)$,
where $\bar G$ is a finite $p$-group, $n\geq2$,  and $\rho\colon G\to \bar G$ is a pro-$p$ group homomorphism.
\end{abstract}


\maketitle

\tableofcontents

{\let\thefootnote\relax\footnote{{This work was supported by the Israel Science Foundation (grant No.\ 569/21).}}}

\section{Introduction}
Let $p$ be a fixed prime number.
Given a profinite group $G$ and an integer $n\geq0$, we write $H^n(G,\dbF_p)$ for the $n$th profinite cohomology group of $G$ with respect to its trivial action on $\dbF_p$, and
let $H^\bullet(G,\dbF_p)=\bigoplus_{n\geq0}H^n(G,\dbF_p)$ be the cohomology (graded) ring with the cup product.
A \emph{symbol} in $H^n(G,\dbF_p)$, where $n\geq1$, is a cup product $\chi_1\cup\cdots\cup\chi_n$ of $n$ elements $\chi_1\nek\chi_n$ of $H^1(G_F,\dbF_p)$.
The \emph{symbol length} $\syml(\omega)$ of a cohomology class $\omega$ in $H^n(G,\dbF_p)$ is the minimal positive integer $m$ such that $\omega$ can be written as a sum of $m$ symbols in $H^n(G,\dbF_p)$;
if no such integer exists, one sets $\syml(\omega)=\infty$.

Now let $F$ be a field which contains a root of unity of order $p$.
Let $G_F=\Gal(F^{\rm sep}/F)$ be its absolute Galois group.
By the celebrated \emph{Norm Residue Theorem}, proved by Voevodsky and Rost,
$H^\bullet(G_F,\dbF_p)$ is isomorphic as a graded ring to the \emph{mod-$p$ Milnor $K$-ring}  $K^M_\bullet(F)/p$ (\cite{Voevodsky11}, \cite{HaesemeyerWeibel19}).
Recall that the latter graded ring is the tensor algebra $\mathrm{Tensor}(F^\times)$ of the multiplicative group $F^\times$ of $F$ divided by the (homogeneous) two-sided ideal generated by all elements $a\tensor (1-a)$ and $a^p$, where $0,1\neq a\in F$.
Consequently, $H^\bullet(G_F,\dbF_p)$ is generated by its degree $1$ elements, that is, $\syml(\omega)<\infty$ for all $\omega$.
General bounds on $\syml(\omega)$ can serve as a measure for the arithmetical complexity of $F$, as well as of specific cohomological constructions (see e.g., \cite{Krashen16}, \cite{Matzri16}).
However, the proof of the Norm Residue Theorem is not constructive, and does not yield such explicit bounds.

In the current paper we consider the symbol length of classes of elements of $H^n(G_F,\dbF_p)$ which are pullbacks of a given cohomology class in a quotient of $G_F$.
More specifically, let $\bar G$ be a fixed pro-$p$ group, and let $\bar\omega$ be a fixed element of $H^n(\bar G,\dbF_p)$.
Let $\rho$ range over all profinite group homomorphisms $G_F\to \bar G$, and write $\rho^*\colon H^n(\bar G,\dbF_p)\to H^n(G_F,\dbF_p)$ for the induced (pullback) homomorphism.
We are interested in the symbol lengths $\syml(\rho^*(\bar\omega))$.
As is explained below, such symbol lengths arise naturally in the study of \emph{Massey products} in $H^2(G_F,\dbF_p)$.

The following conjecture was proposed by Leonid Positselski in a private communication to the author:

\begin{conj}
\label{main conjecture}
Given a pro-$p$ group $\bar G$, a positive integer $n$, and $\bar\omega\in H^n(\bar G,\dbF_p)$, there is a non-negative integer $M=M(\bar G,n,\bar\omega)$ such that for every field $F$ which contains a root of unity of order $p$, and for every profinite group homomorphism $\rho\colon G_F\to \bar G$, one has $\syml(\rho^*(\bar\omega))\leq M$.
\end{conj}

In this conjecture we may replace $G_F$ by its maximal pro-$p$ quotient $G_F(p)$.
Thus $G_F(p)=\Gal(F(p)/F)$, where $F(p)$ is the maximal pro-$p$ Galois extension of $F$.
Indeed, as $H^1(G_{F(p)},\dbF_p)$ is trivial,  the Norm Residue Theorem implies that $H^n(G_{F(p)},\dbF_p)=0$ for every $n\geq1$,
and a standard spectral sequence argument shows that $H^\bullet(G_F(p),\dbF_p)\isom H^\bullet(G_F,\dbF_p)$ via inflation.
In what follows, we will use this alternative pro-$p$ Galois setting.

We note that Conjecture \ref{main conjecture} for $G_F(p)=\bar G$ and $\rho=\id$ implies that every class in $H^n(G_F(p),\dbF_p)$ has a finite symbol length.
Thus it implies the above consequence of the Norm Residue Theorem.

In this paper we prove Conjecture \ref{main conjecture} for an important class of maximal pro-$p$ groups $G_F(p)$, namely, the pro-$p$ groups of \emph{elementary type}.
These are finitely generated pro-$p$ groups which are iteratively constructed from a certain  list of standard \emph{building blocks} using two natural operations: free pro-$p$ products, and certain semi-direct products of the form $\dbZ_p\rtimes\bar G$, called \emph{extensions}.
The iterative construction (recalled in \S\ref{section on cyclotomic pro-p pairs} and \S\ref{section on elementary type pairs}) is in fact carried out in the category of \emph{cyclotomic pro-$p$ pairs} $\calG=(G,\theta)$, which are pro-$p$ groups $G$ with a pro-$p$ group homomorphism  $\theta\colon G\to1+p\dbZ_p$, where $1+p\dbZ_p$ is the group of principal units in the $p$-adic ring $\dbZ_p$.
The pro-$p$ groups of elementary type are then the underlying groups $G$ of the pairs $\calG$ obtained using this process.
Their importance stems from the following

\begin{conj}[The Pro-$p$ Elementary Type Conjecture]
\label{etc}
The finitely generated Galois groups $G_F(p)$, with $F$ as above, are the pro-$p$ groups of elementary type.
\end{conj}

Every pro-$p$ group of elementary type is realizable as a finitely generated maximal pro-$p$ Galois group $G_F(p)$, for some field $F$ containing a root of unity of order $p$ \cite{Efrat98}*{Remark 3.4}, so Conjecture \ref{etc} is actually that the converse also holds.

Our first main result proves Conjecture \ref{main conjecture} for this class of groups:

\begin{thmA}
Given a pro-$p$ group $\bar G$, a positive integer $n$, and a cohomology element $\bar\omega\in H^n(\bar G,\dbF_p)$,  there is a non-negative integer $M=M(\bar G,n,\bar\omega)$ such that for every pro-$p$ group $G$ of elementary type, and every profinite group homomorphism $\rho\colon G\to \bar G$, one has $\syml(\rho^*(\bar\omega))\leq M$.
\end{thmA}

For $\bar G$ finite, such a uniform bound $M$ can be given explicitly -- see Theorem \ref{uniform bound theorem}.

\medskip

The Elementary Type Conjecture originated from a similar conjecture of Marshall in the context of quadratic forms (\cite{Marshall80}, \cite{Marshall04}).
The Galois-theoretic connections in the pro-$2$ case were studied by Jacob and Ware in \cite{JacobWare89}.
The full conjecture in the pro-$p$ Galois context was stated in \cite{Efrat95} (see also \cite{Efrat98}, \cite{Efrat97}).
Among the classes of fields for which this conjecture is known are:
\begin{enumerate}
\item[(1)]
$F$ is a global field (\cite{Efrat97}, \cite{Efrat99a});
\item[(2)]
$F$ is an extension of transcendence degree $\leq1$ of a local field (\cite{Efrat00}, \cite{KrullNeukirch71});
\item[(3)]
$F$ is an extension of transcendence degree $\leq1$ over a pseudo algebraically closed field \cite{Efrat01};
\item[(4)]
$p=2$ and $F$ is a Pythagorean field, i.e., $F^2=F^2+F^2$ (see  \cite{Jacob81},  \cite{Minac86});
\item[(5)]
$p=2$ and $G_F(2)$ is generated as a pro-$2$ group by $\leq5$ elements \cite{CarsonMarshall82}.
\end{enumerate}

Hence Conjecture \ref{main conjecture} holds for these classes of fields if $G_F(p)$ is finitely generated.

\medskip

Since the Norm Residue Theorem gives a complete elementary description of the graded ring $H^\bullet(G_F,\dbF_p)=H^\bullet(G_F(p),\dbF_p)$, for $F$ containing a root of unity of order $p$, much of the research in recent years on the cohomology of $G_F$ has focused on \emph{external} operations on
$H^\bullet(G_F,\dbF_p)$, which go beyond its ring structure, and mainly on \emph{Massey products}, which we now describe.

While Massey products are defined in the general homological context of differential graded algebras (also in higher degrees), in the special context of profinite group cohomology and the 2nd cohomology group, they have the following equivalent definition, due to Dwyer \cite{Dwyer75} (in the discrete case):
For an integer $m\geq2$,  let $\dbU_m(\dbF_p)$ be the group of all unipotent upper-triangular $(m+1)\times(m+1)$-matrices over $\dbF_p$.
Let $\overline{\dbU}_m(\dbF_p)$ be the quotient of $\dbU_m(\dbF_p)$ by its subgroup consisting of all matrices which are $0$ except for the main diagonal and the upper-right entry.
There is a central extension of finite $p$-groups
\[
1\to \dbF_p\to\dbU_m(\dbF_p)\to\overline{\dbU}_m(\dbF_p)\to1.
\]
It corresponds to a cohomology class $\bar\omega_m\in H^2(\overline\dbU_m(\dbF_p),\dbF_p)$.
Now let $G$ be a profinite group and let $\chi_1\nek\chi_m\in H^1(G,\dbF_p)=\Hom_{\rm cont}(G,\dbF_p)$.
The \emph{$m$-fold Massey product} $\langle\chi_1\nek\chi_m\rangle$ of $\chi_1\nek\chi_m$  is the set of all pullbacks $\rho^*(\bar\omega_m)\in H^2(G,\dbF_p)$,
where $\rho\colon G\to\overline{\dbU}_m(\dbF_p)$ is a profinite group homomorphism whose projections to the super-diagonal are $\chi_1\nek\chi_m$, i.e., $\proj_{i,i+1}\circ\rho=\chi_i$, $i=1,2\nek m$.
When $m=2$ one has $\langle\chi_1,\chi_2\rangle=\{\chi_1\cup\chi_2\}$, so in this sense, Massey products extend the usual ring structure of $H^\bullet(G,\dbF_p)$.

In \cite{MinacTan15} Min\'a\v c and T\^an proposed the following
\begin{conj}
[The $m$-Massey Vanishing Conjecture]
\label{the Massey vanishing conjecture}
Let $F$ be a field which contains a root of unity of order $p$, let $m\geq3$, and let $\chi_1\nek\chi_m\in H^1(G_F(p),\dbF_p)$.
If the $m$-fold Massey product $\langle\chi_1\nek\chi_m\rangle$ is nonempty, then it contains $0$.
\end{conj}

Conjecture \ref{the Massey vanishing conjecture} was proved in particular in the following cases:
\begin{enumerate}
\item[(1)]
$m=3$, $p=2$, and $F$ is a local or a global field (Hopkins and Wickelgren \cite{HopkinsWickelgren15});
\item[(2)]
$m=3$ and $p$ and $F$ are arbitrary (Matrzi \cite{Matzri14}, followed by Efrat--Matzri \cite{EfratMatzri17} and  Min\'a\v c--T\^an \cite{MinacTan16}, and more recently, by Hay et al.\ \cite{HayLamSharifiWangWake23});
\item[(3)]
$m=4$, $p=2$, and $F$ arbitrary (Merkurjev and Scavia \cite{MerkurjevScavia23});
\item[(4)]
$m\geq3$ and $F$ is a number field (Harpaz and Wittenberg \cite{HarpazWittenberg23}).
\end{enumerate}

Trivially, for every homomorphism $\rho\colon G_F(p)\to\overline{\dbU}_2(\dbF_p)$ one has $\syml(\rho^*(\bar\omega_2))\leq1$.
The (proven) case $m=3$ of the Vanishing Conjecture means that, for $\chi_1,\chi_2,\chi_3\in H^1(G_F(p),\dbF_p)$, every nonempty Massey product $\langle\chi_1,\chi_2,\chi_3\rangle$ coincides with its \emph{indeterminacy},
i.e., it consists of all  sums of the form $\chi_1\cup\psi+\psi'\cup\chi_3$ with $\psi,\psi'\in H^1(G_F(p),\dbF_p)$ \cite{EfratMatzri17}*{\S1}.
Therefore the pullbacks $\rho^*(\bar\omega_3)$ satisfy $\syml(\rho^*(\bar\omega_3))\leq2$.
It is a special case of Conjecture \ref{main conjecture} that for every $m$ there should be a \emph{uniform bound}, depending only on $m$, on the symbol length of all $m$-fold Massey product elements $\rho^*(\bar\omega_m)$ in $H^2(G_F(p),\dbF_p)$.
Thus, this special case  -- which in fact motivated the current work -- partly generalizes the proven case $m=3$ of the Massey vanishing conjecture.
On the other hand, for $m\geq4$ there appears to be no direct connection between Conjecture \ref{main conjecture}, which is the focus of this paper, and the Massey vanishing conjecture. 

As an application of Theorem A, we prove in \S\ref{section on Massey products}:

\begin{thmB}
Let $m\geq2$, and let $G$ be a pro-$p$ group of elementary type.
Then the elements of every $m$-fold Massey product in $H^2(G,\dbF_p)$ have symbol length at most $\lfloor m^2/4\rfloor+m$.
\end{thmB}

In the special case of pro-$p$ groups of elementary type, Conjecture \ref{the Massey vanishing conjecture} was proved for every $m\geq3$ by Quadrelli \cite{Quadrelli24} (when $p=2$ one also needs to assume that $\sqrt{-1}\in F$).
Moreover, he proved for this class of groups a stronger fact,
where the non-emptiness of $\langle\chi_1\nek\chi_m\rangle$ is replaced by the weaker assumption that $\chi_i\cup\chi_{i+1}=0$, $i=1,2\nek m-1$.
For $m=3$ the two versions are equivalent \cite{EfratMatzri17}*{\S1}, however it was shown by Harpaz and Wittenberg that already for $m=4$ and $p=2$, the stronger version does not hold for $G_{\dbQ}(2)$ (see \cite{GuillotMinacTopaz18}*{Appendix, Example A.15 and p.\ 1950}).
This was recently extended by Merkurjev and Scavia for more general fields \cite{MerkurjevScavia22}*{Th.\ 1.6}.
We refer to \cite{HarpazWittenberg23} for more references and information on the history of Conjecture \ref{the Massey vanishing conjecture}.

It should be mentioned that several other major conjectures on the structure of absolute and maximal pro-$p$ Galois groups were recently proved for the class of pro-$p$ groups of elementary type:
In \cite{MinacPasiniQuadrelliTan21} the authors verify for a group $G$ in this class a conjecture of Positselski asserting that the cohomology algebra $H^\bullet(G,\dbF_p)$ is Koszul \cite{Positselski14}*{\S0.1} (see also \cite{Positselski05}, \cite{PositselskiVishik95}).
Various stronger variants of Koszulity  are similarly proved for elementary type groups in \cite{MinacPasiniQuadrelliTan22}.
Another conjecture, the \emph{Bogomolov--Positselski freeness conjecture} on the maximal pro-$p$ Galois group of the maximal $p$-radical extension of a field \cite{Positselski05}*{Conjecture 1.2}, was proved for elementary type pro-$p$ groups by Quadrelli and Weigel \cite{QuadrelliWeigel22}.

\bigskip

We conclude the Introduction by sketching the proofs of the main theorems.

Let $\bar G$, $G$, $\rho\colon G\to\bar G$ and $\bar\omega$ be as in Theorem A.
Thus there is a cyclotomic pro-$p$ pair $\calG=(G,\theta)$ of elementary type.
By a limit argument, we may assume that $\bar G$ is a finite $p$-group.
For simplicity, we first assume that $n=2$.
The basic observation behind the proof is that the maximal symbol length of elements of $H^2(G,\dbF_p)$ does not grow under free products of pro-$p$ pairs, and grows by $1$ under the operation of extension (see \S\ref{section on the symbol length under elementary operations}).
We therefore consider maximal chains of extensions
\[
Z_r\rtimes(Z_{r-1}\rtimes(\cdots (Z_1\rtimes\calB)\cdots))
\]
which are involved in the inductive construction of $\calG$, with $Z_1\isom\cdots\isom Z_r\isom\dbZ_p$ and $\calB$ a building block.
Then $V^j=Z_{j+1}\times Z_{j+2}\times\cdots\times Z_r$ is a free pro-$p$ abelian closed subgroup of $G$ for each $0\leq j\leq r$ (where  $V^r=\{1\}$).

Let $l(\bar G)$ be the maximal integer $l\geq0$ such that $\bar G$ contains a proper chain of $l+1$ abelian subgroups.
If $r>l(\bar G)$, then $\rho(V^k)=\rho(V^{k+1})$ for some $0\leq k\leq r-1$.
There is a pro-$p$ group automorphism $\alp$ of $V^0$ such that $\alp(V^j)=V^j$, $j=1,2\nek r$, and $\rho\circ\alp$ is trivial on $V^k$ (Proposition \ref{normalization lemma}).
In \S\ref{section on automorphisms} we show that $\alp$ is induced by an automorphism of $\calG$, so we may assume without loss of generality that $\rho(V^k)=\{1\}$.
As shown in \S\ref{section on factoring of epimorphisms}, $\rho$ then factors via a proper quotient $\calG'$ of $\calG$ which is also of elementary type.
By induction, we obtain such a quotient for which the maximal length $r$ of a chain of extensions as above is at most $l(\bar G)$.
Now symbol lengths of pullbacks $\rho^*\bar\omega$ can only increase when passing from $\calG$ to $\calG'$, and
furthermore, the symbol lengths at the level of the building blocks $\calB$ are at most $1$.
Altogether this gives the bound
\[
\syml(\rho^*(\bar\omega))\leq r+1\leq l(\bar G)+1.
\]

When $n>2$ the symbol lengths can grow by more than $1$ under extensions.
An analysis due to Wadsworth \cite{Wadsworth83} allows us to bound the possible growth (\S\ref{section on the symbol length under elementary operations}), and amend the above upper bound accordingly.

For the proof of Theorem B, we use known upper bounds, due to Goozeff \cite{Goozeff70} and Barry \cite{Barry 79}, on the order of abelian subgroups of $\dbU_m(\dbF_p)$, to derive similar bounds for $\overline\dbU_m(\dbF_p)$.
Theorem A then yields the required bound on the symbol length of Massey product elements $\rho^*(\bar\omega_m)$ for pro-$p$ groups of elementary type.

\begin{acknowledgment}
\rm
This paper grew out of discussions and correspondences with Leonid Positselski.
I acknowledge with gratitude his substantial contribution, especially in the strategy of the proof of Theorem A.
I thank the referee for the careful reading of this manuscript and the very helpful comments and suggestions.
\end{acknowledgment}

\section{Cyclotomic Pro-$p$ Pairs}
\label{section on cyclotomic pro-p pairs}
Throughout this paper we fix a prime number $p$.
Let $\dbZ_p^\times$ be the multiplicative group of the ring $\dbZ_p$ of $p$-adic integers, and let $1+p\dbZ_p$ be its subgroup of principal units.
 A \emph{cyclotomic pro-$p$ pair} \cite{Efrat98} (or simply, \emph{``a pro-$p$ pair''}) $\calG=(G,\theta)$ consists of a pro-$p$ group $G$
and a pro-$p$ group homomorphism $\theta\colon G\to 1+p\dbZ_p$.
A \emph{morphism} $\phi\colon \calG_1=(G_1,\theta_1)\to\calG_2=(G_2,\theta_2)$ of pro-$p$ pairs is a pro-$p$ group homomorphism $\phi\colon G_1\to G_2$ such that $\theta_1=\theta_2\circ\phi$.
A \emph{pro-$p$ subpair} of $\calG$ is a pro-$p$ pair $(G_0,\theta|_{G_0})$, where $G_0$ is a closed subgroup of $G$.

The \emph{free product} of pro-$p$ pairs $\calG_1=(G_1,\theta_1)$ and $\calG_2=(G_2,\theta_2)$ is the pro-$p$ pair
\[
\calG_1*\calG_2=(G_1*_pG_2,\>\theta_1*_p\theta_2),
\]
where $G_1*_pG_2$ is the free product taken in the category of pro-$p$ groups, and $\theta_1*_p\theta_2\colon G_1*_pG_2\to 1+p\dbZ_p$ is the unique homomorphism of pro-$p$ groups which extends $\theta_1$ and $\theta_2$, as provided by the universal property of $G_1*_pG_2$.
This construction has the following universal property:
For any cyclotomic pro-$p$ pair $\calK$, any two morphisms $\calG_1\to\calK$ and $\calG_2\to\calK$ uniquely extend to a morphism
$\calG_1*\calG_2\to\calK$.

Next let $Z$ be a free abelian pro-$p$ group.
Thus $Z\cong\dbZ_p^m$ for some cardinal number $m$, so we have a multiplication map $\dbZ_p^*\times Z\to Z$.
We will henceforth use multiplicative notation for $Z$.
For  a cyclotomic pro-$p$ pair $\bar\calG=(\bar G,\bar\theta)$, let  $Z\rtimes \bar G$ be the semi-direct product with respect to $\bar\theta$;
That is, for $\bar g\in\bar G$ and $z\in Z$ we set
\[
\bar gz\bar g^{-1}=z^{\bar\theta(\bar g)}.
\]
The \emph{extension} of $\bar\calG$ by $Z$ is the cyclotomic pro-$p$ pair
\[
Z\rtimes\bar\calG=(Z\rtimes \bar G,\>\theta),
\]
where $\theta\colon Z\rtimes \bar G\to 1+p\dbZ_p$ is the composition of the projection $\pi\colon Z\rtimes \bar G\to\bar G$ with $\bar\theta\colon\bar G\to 1+p\dbZ_p$.
Then $\pi\colon\calG\to\bar\calG$ is an epimorphism of pro-$p$ pairs.
The pro-$p$ pair $\bar\calG$ embeds as a subpair of $Z\rtimes\bar\calG$ in a natural way.
If $Z'$ is another free abelian pro-$p$ group, then we identify $(Z\times Z')\rtimes\bar\calG=Z\rtimes(Z'\rtimes\bar\calG)$.

This construction has the following universal property:
Let $\calK=(K,\kappa)$ be a pro-$p$ pair.
Suppose further that  $\psi\colon Z\to \Ker(\kappa)$ is a pro-$p$ group homomorphism and $\bar\eta\colon\bar\calG\to\calK$ is a morphism,
such that for every $\bar g\in\bar G$ and $z\in Z$ one has  $\bar\eta(\bar g)\psi(z)\bar\eta(\bar g)^{-1}=\psi(z)^{\bar\theta(\bar g)}$.
Then there is a unique morphism $\delta\colon Z\rtimes\bar\calG\to\calK$ such that $\delta|_Z=\psi$ and $\delta|_{\bar\calG}=\bar\eta$.

In particular, any morphism $\phi\colon\bar\calG\to\bar\calG'$ of pro-$p$ pairs and any homomorphism $\psi\colon Z\to Z'$ of free abelian pro-$p$ groups, induce a morphism $\psi\rtimes\phi\colon Z\rtimes\bar\calG\to Z'\rtimes\bar\calG'$.

\begin{rem}
\label{extended semidirect product relation}
\rm
The above defining relation of the semi-direct product can be slightly extended as follows:
Write $\calG=Z\rtimes\bar\calG=(G,\theta)$.
For every $g\in G=Z\rtimes\bar G$ and $z\in Z$ we have
\[
gzg\inv=z^{\theta(g)}.
\]
Indeed, write $g=\bar gz'$ with $\bar g\in\bar G$ and $z'\in Z\subseteq\Ker(\theta)$.
Then $z$ and $z'$ commute, and $\theta(g)=\bar\theta(\bar g)$.
Therefore
$gzg\inv=\bar gz\bar g\inv=z^{\bar\theta(\bar g)}=z^{\theta(g)}$.
\end{rem}

\begin{lem}
\label{morphisms into extensions}
Let $\phi\colon\bar\calG'\to\bar\calG$ be a morphism of pro-$p$ pairs, and let $\psi\colon Z'\to Z$ be a homomorphism of free abelian pro-$p$ groups.
Let $\pi\colon Z\rtimes\bar\calG\to\bar\calG$ and $\pi'\colon Z'\rtimes\bar\calG'\to\bar\calG'$ be the projection morphisms, and suppose that $\bar\eta\colon\bar\calG'\to Z\rtimes\bar\calG$ is a morphism such that $\pi\circ\bar\eta=\phi$.
Then:
\begin{enumerate}
\item[(a)]
There is a unique morphism $\eta\colon Z'\rtimes\bar\calG'\to Z\rtimes\bar\calG$ which extends $\psi$ and $\bar\eta$.
\item[(b)]
One has $\pi\circ\eta=\phi\circ\pi'$, i.e., the outer square in the following diagram commutes:
\[
\xymatrix{
Z\rtimes\bar\calG\ar[r]^{\ \ \pi}&\bar\calG\\
Z'\rtimes\bar\calG'\ar[u]^{\eta}\ar[r]^{\ \ \pi'}&\bar\calG'\ar[ul]_{\bar\eta}\ar[u]_{\phi}.\\
}
\]
\item[(c)]
If in addition $\phi$ and $\psi$ are injective (resp., surjective), then $\eta$ is injective (resp., surjective).
\end{enumerate}
\end{lem}
\begin{proof}
Write $\bar\calG=(\bar G,\bar\theta)$ and $\bar\calG'=(\bar G',\bar\theta')$.
Then $Z\rtimes\bar\calG=(Z\rtimes\bar G,\theta)$, where $\theta$ is induced by $\bar\theta$.

(a) \quad
Consider $z'\in Z'$ and $\bar g'\in\bar G'$.
By Remark \ref{extended semidirect product relation}, and since $\bar\eta$ is a morphism,
\[
\bar\eta(\bar g')\psi(z')\bar\eta(\bar g')\inv=\psi(z')^{\theta(\bar\eta(\bar g'))}=\psi(z')^{\bar\theta'(\bar g')}.
\]
The existence and uniqueness of $\eta$ follow from the universal property of $Z'\rtimes\bar\calG'$.

\medskip

(b)\quad
We observe that $\pi\circ\eta$ and $\phi\circ\pi'$ are trivial on $Z'$ and, by assumption, coincide on $\bar G'$.
Hence they coincide on $Z'\rtimes\bar\calG'$.

\medskip

(c) \quad
Suppose that $\phi$ and $\psi$ are injective, let $z'\in Z'$ and $\bar g'\in\bar G'$, and assume that $\eta(z'\bar g')=1$.
By (b), $\phi(\bar g')=\phi(\pi'(z'\bar g'))=\pi(\eta(z'\bar g'))=1$, so $\bar g'=1$.
Hence $\psi(z')=\eta(z')=1$, and therefore $z'=1$.

Finally suppose that $\phi$ and $\psi$ are surjective.
By (b),  $\pi\circ\eta=\phi\circ\pi'$ is also surjective.
Further, $\Ker(\pi)=Z=\Img(\psi)\subseteq \Img(\eta)$.
Hence, $\eta$ is surjective.
\end{proof}

\medskip

We now describe several important examples of cyclotomic pro-$p$ pairs.

\begin{exam}
\label{pro-p pair of a field}
\rm
Let $F$ be a field which contains a root of unity of order $p$.
For every $i\ge1$ let $\mu_{p^i}$ be the group of all $p^i$th roots of unity in the algebraic closure of  $F$;
In fact, $\mu_{p^i}$ is contained in the maximal pro-$p$ Galois extension $F(p)$ of $F$.
Let $\mu_{p^\infty}=\bigcup_{i\geq1}\mu_{p^i}$.
As before, let  $G_F(p)=\Gal(F(p)/F)$.
The \emph{pro-$p$ cyclotomic character} $\theta_{F,p}$ of $F$ is the composition of the restriction homomorphism $G_F(p)\to\Aut(\mu_{p^\infty}/\mu_p)$ with the canonical isomorphism $\Aut(\mu_{p^\infty}/\mu_p)\cong1+p\dbZ_p\,(\leq\dbZ_p^\times)$.
Then the \emph{Galois pro-$p$ pair of $F$} is
\[
\calG_F(p)=\bigl(G_F(p),\theta_{F,p}\bigr).
\]
\end{exam}

\begin{exam}
\label{the standard building blocks}
\rm
The following cyclotomic pro-$p$ pairs will be called the \emph{standard building blocks} in this category:
\begin{enumerate}
\item[(1)]
The \emph{trivial} pro-$p$ pair $(1,1)$;
\item[(2)]
All pro-$p$ pairs $(\dbZ_p,\theta)$, where $\theta\colon\dbZ_p\to1+p\dbZ_p$ is a pro-$p$ group homomorphism;
\item[(3)]
The Galois pro-$p$ pairs $\calG_F(p)$, where $F$ is a finite extension of $\dbQ_p$ containing a root of unity of order $p$;
\item[(4)]
When $p=2$, the pro-$2$ pair $\calE=(\dbZ/2,\theta)$, where $\theta$ is the unique nontrivial homomorphism $\dbZ/2\to\{\pm1\}\subset 1+2\dbZ_2$.
Note that $\calE=\calG_{\dbR}$.
\end{enumerate}
\end{exam}

\begin{rem}
\label{Demushkin-groups-remark}
\rm
While this will not be needed in the sequel, we remark that the structure of the Galois pro-$p$ pairs in (3) can be described group-theoretically.
For this, we recall that a pro-$p$ group $G$ is a \emph{Demushkin group}
if $H^1(G,\dbF_p)$ is finite, $\dim_{\dbF_p}H^2(G,\dbF_p)=1$, and the cup product
\[
\cup\colon H^1(G,\dbF_p)\times H^1(G,\dbF_p)\to H^2(G,\dbF_p)
\]
is non-degenerate.
The only finite Demushkin group is $\dbZ/2$, when $p=2$ \cite{NeukirchSchmidtWingberg}*{Prop.\ 3.9.10}.
When $G$ is infinite, these properties further imply that $H^n(G,\dbF_p)=0$ for all $n\geq3$ \cite{NeukirchSchmidtWingberg}*{Th.\ 3.7.2}.
Every Demushkin group $G$ is equipped with a canonical continuous homomorphism $\chi\colon G\to1+p\dbZ_p$, characterized by the property that the natural map $H^1(G,I/p^iI)\to H^1(G,I/pI)$ is surjective for all $i\geq 1$, where $I=\dbZ_p(\chi)$ is the discrete $G$-module $\dbZ_p$ with action given by $\chi$ \cite{Labute67}*{Th.\ 4}.
Works by Demushkin \cite{Demushkin61}, Serre \cite{Serre63}, and Labute \cite{Labute67} (see also \cite{NeukirchSchmidtWingberg}*{\S3.9}) give an explicit presentation of any Demushkin pro-$p$ group $G$ in terms of (pro-$p$) generators and a single relation, as well as the values of $\chi$ on these generators.

Now let $F$ be a field containing a root of unity of order $p$, and let $G=G_F(p)$.
Then the map $H^1(G,I/p^iI)\to H^1(G,I/pI)$ is surjective when we take $\chi=\chi_{F,p}$.
Indeed, then $I/p^iI=\mu_{p^i}$ and $I/pI=\mu_p$, so by Kummer theory, this map may be identified with the natural projection $F^\times/(F^\times)^{p^i}\to F^\times/(F^\times)^p$, which is obviously surjective.
If, moreover, $F$ is a finite extension of $\dbQ_p$, then local class field theory implies that $G=G_F(p)$ is a Demushkin group \cite{NeukirchSchmidtWingberg}*{Th.\ 7.5.11}.
Therefore the above works by Demushkin, Serre, and Labute provide a complete group-theoretic description of the Galois pro-$p$ pairs $\calG_F(p)$ in case (3).

It is currently an open problem whether there exists any other example of a field $F$ containing a root of unity of order $p$ such that $G_F(p)$ is a Demushkin group.
\end{rem}

\section{Elementary Type Pro-$p$ Pairs}
\label{section on elementary type pairs}
In the rest of the paper, we fix a nonempty class $\grB$ of cyclotomic pro-$p$ pairs, called the \emph{(generalized) building blocks}.
The class $\mathrm{ET}(\grB)$ of  \emph{elementary type pro-$p$ pairs over $\grB$} will consist of all cyclotomic pro-$p$ pairs which can be constructed in finitely many steps from the pairs in $\grB$ using free products and extensions.

We will often need to keep track of the way in which a particular  elementary type
pro-$p$ pair is constructed from building blocks.
Thus we will speak about \emph{elementary type constructions} $\cc$ over $\grB$, and to each such construction we associate its pro-$p$ pair $\calG(\cc)$.
These constructions will be formal expressions, defined inductively as follows:
\begin{enumerate}
\item[(i)]
Every $\calB\in\grB$ is an elementary type construction, and we set $\calG(\calB)=\calB$;
\item[(ii)]
If $\cc_1,\cc_2$ are elementary type constructions different from $(1,1)$, then we have an elementary type construction $\cc_1\diamond\cc_2$ with $\calG(\cc_1\diamond\cc_2)=\calG(\cc_1)*\calG(\cc_2)$;
\item[(iii)]
If $\bar\cc$ is an elementary type construction, then we have an elementary type construction $\langle\bar \cc\rangle$ with $\calG(\langle\bar \cc\rangle)=\dbZ_p\rtimes\calG(\bar\cc)$.
\end{enumerate}
Thus technically, the elementary type constructions are certain words in the alphabet $\grB\cup\{\diamond,\langle,\rangle,(,)\}$.
The \emph{elementary type pro-$p$ pairs} over $\grB$ will be the pro-$p$ pairs $\calG(\cc)$, where $\cc$ is an elementary type construction $\cc$ over $\grB$.
The \emph{pro-$p$ groups of elementary type over $\grB$} are the underlying groups of elementary type pro-$p$ pairs.

\begin{exam}
\label{c to G is not injective}
\rm
Let $\calB$ be a nontrivial building block in $\grB$.
The elementary type constructions
$\cc=\calB\diamond\calB$ and $\cc'=\langle\calB\rangle$
have associated pro-$p$ pairs
\[
\calG(\cc)=\calB*\calB, \quad \calG(\cc')=\dbZ_p\rtimes\calB.
\]
For $p=2$ and $\calB=\calE$ as in Example \ref{the standard building blocks}(4), $\calG(\cc),\calG(\cc')$ are isomorphic pro-$2$ pairs.
\end{exam}

\begin{exam}
\label{standard elementary type pairs example}
\rm
Let $\grB$ be the class of all standard building blocks, as in Example \ref{the standard building blocks}.
Then $\mathrm{ET}(\grB)$ is the class of \emph{standard} elementary type pro-$p$ pairs.
By induction, the elementary type pro-$p$ groups in this case are finitely generated.
Every pro-$p$ pair in $\mathrm{ET}(\grB)$ is realizable as $\calG_F(p)$ for some field $F$ of characteristic $0$ which contains a root of unity of order $p$ \cite{Efrat98}*{Remark 3.4}.
As discussed in the Introduction, the \emph{Elementary Type Conjecture} predicts that, conversely, if $F$ is a field containing a root of unity of order $p$ and $G_F(p)$ finitely generated as a pro-$p$ group, then $\calG_F(p)\in\mathrm{ET}(\grB)$.
\end{exam}

\begin{exam}
\label{absolute elementary type pairs example}
\rm
Let $\grB$ be the class of all standard building blocks of the forms (1), (2) and (4) of Example \ref{the standard building blocks}.
We call $\mathrm{ET}(\grB)$ the class of \emph{absolute} elementary type pro-$p$ pairs.
Every such pro-$p$ pair is realizable as $\calG_F(p)$ for some field $F$ of characteristic $0$ such that the absolute Galois group $G_F$ is a finitely generated pro-$p$ group (whence $F$ contains a root of unity of order $p$) \cite{Efrat98}*{Remark 3.4}.
A similar conjecture predicts that, conversely,  if $F$ is a field such that $G_F$ is a finitely generated pro-$p$ group, then $\calG_F(p)\in\mathrm{ET}(\grB)$.
\end{exam}

\begin{exam}
\label{elementary type pairs with all roots of unity example}
\rm
Let $\grB=\{(1,1)\}$.
Then every pro-$p$ pair in $\mathrm{ET}(\grB)$ is realizable as $\calG_F(p)$ for some field $F$ of characteristic $0$ which contains $\mu_{p^\infty}$ with $G_F(p)$ a finitely generated pro-$p$ group.
Conversely, it is a special case of the Elementary type Conjecture that, if $F$ is a field of characteristic $\neq p$ which contains $\mu_{p^\infty}$ and $G_F(p)$ is finitely generated, then $\calG_F(p)\in \mathrm{ET}(\grB)$.
\end{exam}

\begin{exam}
\label{pythagorean elementary type pairs example}
\rm
Assume that $p=2$ and let $\grB=\{\calE\}$, where $\calE$ is as in (4) of Example \ref{the standard building blocks}.
By results of Marshall \cite{Marshall80}, Jacob \cite{Jacob81}, Craven \cite{Craven78}, and Min\'a\v c \cite{Minac86}, $\mathrm{ET}(\grB)$ is the class of all Galois pro-$2$ pairs $\calG_F(p)$,
where $F$ is a real-Pythagorean field, i.e., $-1\not\in F^2=F^2+F^2$, and $G_F(2)$ is a finitely generated pro-$2$ group (See \cite{EfratHaran94}*{\S3}).
\end{exam}

There is a natural partial order $\leq$ on the set of all elementary type constructions over $\grB$, where $\dd\leq\cc$ means that $\dd$ is an \emph{subconstruction} of $\cc$.
More formally, we define this relation by induction on the structure of $\cc$ as follows:
\begin{enumerate}
\item[(i)]
The only subconstruction of $\cc=\calB$, where $\calB\in\grB$, is $\calB$ itself;
\item[(ii)]
The subconstructions of $\cc=\cc_1\diamond\cc_2$ are the subconstructions $\cc'_1$, $\cc'_2$ of $\cc_1,\cc_2$, respectively, as well as $\cc'_1\diamond\cc'_2$, if $\cc'_1,\cc'_2\neq(1,1)$;
\item[(iii)]
The subconstructions of  $\cc=\langle\bar\cc\rangle$ are the elementary type constructions $\bar\cc'$ and $\langle\bar\cc'\rangle$, where $\bar\cc'$ ranges over all subconstructions of $\bar\cc$.
\end{enumerate}

Note that a given elementary type construction $\cc$ has only finitely many subconstructions.
This allows inductive arguments on the structure of $\cc$.

\medskip

\begin{rem}
\label{embedding of subconstruction}
\rm
Given  elementary type constructions $\dd\leq\cc$ over $\grB$,
there are a natural monomorphism of pro-$p$ pairs $\iota_{\dd,\cc}\colon \calG(\dd)\hookrightarrow\calG(\cc)$, and a natural epimorphism of pro-$p$ pairs  $\pi_{\cc,\dd}\colon \calG(\cc)\hookrightarrow\calG(\dd)$.

Indeed, by the above universal properties, if $\calK_1,\calK_2,\bar\calK$ are subpairs (resp., epimorphic images) of pro-$p$ pairs $\calG_1,\calG_2,\bar \calG$, then $\calK_1*\calK_2$ is a subpair (resp., epimorphic image)  of $\calG_1*\calG_2$, and $\dbZ_p\rtimes\bar\calK$ is a subpair (resp., epimorphic image) of $\dbZ_p\rtimes\bar\calG$, in a natural way.

We observe that $\pi_{\cc,\dd}\circ\iota_{\dd,\cc}=\id_{\calG(\dd)}$.

More formally, the partially ordered set of elementary type constructions over $\grB$ forms a category, with a single morphism $\dd\to\cc$ when $\dd\leq\cc$.
Then $\cc\mapsto\calG(\cc)$ provides a covariant (resp., contravariant) functor from this category to the category of pro-$p$ pairs, where the morphism $\dd\to\cc$ is mapped to $\iota_{\dd,\cc}$  (resp., to $\pi_{\cc,\dd}$).
\end{rem}

\section{Principal Tuples}
\label{section on principal tuples}

We consider tuples $\calA=(Z_1\nek Z_r)$ of copies $Z_1\nek Z_r$ of $\dbZ_p$.
We call $\calA$ a \emph{$\dbZ_p$-tuple}, and refer to the nonnegative integer $r=\rk(\calA)$ as its \emph{rank}.
Thus $A:=Z_1\times\cdots\times Z_r$ is then a free abelian pro-$p$ group of rank $r$.
We use multiplicative notation for the $Z_i$ and $A$.
When $r=0$, we call $\calA$ the \emph{trivial} $\dbZ_p$-tuple.
We set
\[
V^j\calA=Z_{j+1}\times Z_{j+2}\times\cdots\times Z_r, \quad j=0,1\nek r,
\]
where $V^r\calA=\{1\}$.
Thus $V^0\calA=A$.

Given a pro-$p$ group endomorphism $\alp$ of $A$ and $1\leq k\leq r$, let $\bar\alp_k$ be the endomorphism of $Z_1\times\cdots\times Z_k$ given by the composition
\[
Z_1\times\cdots\times Z_k\xrightarrow{\alp|_{Z_1\times\cdots\times Z_k}}
A=Z_1\times\cdots\times Z_r\to A/V^k\calA=Z_1\times\cdots\times Z_k.
\]
Also let $\alp_k$ be the restriction of $\alp$ to $Z_{k+1}\times\cdots\times Z_r$.

We say that a pro-$p$ group automorphism $\alpha$ of $A=V^0\calA$ is an \emph{$\calA$-automorphism} if $\alpha(V^j\calA)=V^j\calA$ for every $1\leq j\leq r$.
The $\calA$-automorphisms form a subgroup of the group of all pro-$p$ group automorphisms of $A$.

The proof of the following fact is straightforward:

\begin{lem}
\label{induced automorphisms}
Let $\calA=(Z_1\nek Z_r)$ be a $\dbZ_p$-tuple, let $\alpha$ be an automorphism of the abelian pro-$p$ group $A=V^0\calA$, and let $1\leq k\leq r$.
Then $\alpha$ is an $\calA$-automorphism if and only if $\bar\alpha_k$ is a $(Z_1\nek Z_k)$-automorphism, and $\alpha_k$ takes $Z_{k+1}\times\cdots\times Z_r$ into itself and is a $(Z_{k+1}\nek Z_r)$-automorphism.
\end{lem}

\begin{prop}
\label{normalization lemma}
Let $\calA=(Z_1\nek Z_r)$ be a $\dbZ_p$-tuple, let $\bar G$ be a pro-$p$ group, and let $\rho_0\colon A=V^0\calA\to \bar G$ be a pro-$p$ group homomorphism.
For $0\leq k\leq r-1$, the following conditions are equivalent:
\begin{enumerate}
\item[(a)]
$\rho_0(V^k\calA)=\rho_0(V^{k+1}\calA)$;
\item[(b)]
There is an $\calA$-automorphism $\alp$ of $A$ such that $\rho_0(\alp(Z_k))=\{1\}$.
\end{enumerate}
\end{prop}
\begin{proof}
(a)$\Rightarrow$(b): \quad
Choose generators $u_1\nek u_r$ of the pro-$p$ groups $Z_1\nek Z_r$, respectively.
By assumption, there exists $v\in V^{k+1}\calA$ with $\rho_0(u_k)=\rho_0(v)$.
Then $w=v\inv u_k$ is an element of $V^k\calA$ whose image generates $V^k\calA/V^{k+1}\calA$ as a pro-$p$ group, and further, $\rho_0(w)=1$.
Therefore there exists $t\in\dbZ_p^\times$ such that $wu_k^{-t}\in V^{k+1}\calA$.
It follows that the elements $u_1\nek u_{k-1},w,u_{k+1}\nek u_r$ form a basis for the $\dbZ_p$-module $A$.
We define the  $\calA$-automorphism $\alp$ of $A$ by $\alp(u_k)=w$ and $\alp(u_j)=u_j$ for $j\neq k$.
Note that $\rho_0(\alp(Z_k))=\langle\rho_0(\alp(u_k))\rangle=\langle\rho_0(w)\rangle=\{1\}$.

\medskip

(b)$\Rightarrow$(a): \quad
Since $V^k\calA=Z_k\times V^{k+1}\calA$ and $\alpha$ is an $\calA$-automorphism of $A$, we have $V^k\calA=\alpha(Z_k)\times V^{k+1}\calA$, and (a) follows from $\rho_0(\alp(Z_k))=\{1\}$.
\end{proof}

Next, to any elementary type construction $\cc$ over $\grB$ with a pro-$p$ pair $\calG(\cc)=(G,\theta)$, we associate a set of \emph{principal tuples}.
These will consist of a $\dbZ_p$-tuple $(Z_1\nek Z_r)$  such that $Z_1\times\cdots\times Z_r$ is a closed subgroup of $G$, together with a building block $\calB\in\grB$ which appears in $\cc$, and which we call the \emph{root} of the principal tuple.
Roughly speaking, such $\dbZ_p$-tuples will arise from a maximal sequence of iterated extensions in the construction of $\calG(\cc)$, starting from $\calB$.
More specifically, the definition of this set is by induction on the structure of $\cc$, as follows:

\medskip

(i) \quad
If $\cc=\calB$, where $\calB\in\grB$ is a building block, then the unique principal tuple of $\calG(\cc)$ is the trivial $\dbZ_p$-tuple, with the root $\calB$.

(ii) \quad
Suppose that $\cc=\cc_1\diamond\cc_2$ for nontrivial elementary type constructions $\cc_1,\cc_2$ over $\grB$.
Then the principal tuples of $\calG(\cc)=\calG(\cc_1)*\calG(\cc_2)$ will be the principal tuples of either $\calG(\cc_1)$ or $\calG(\cc_2)$, with the same root.

(iii) \quad
Suppose that $\cc=\langle\bar\cc\rangle$ for an elementary type construction $\bar\cc$ over $\grB$.
Then $\calG(\cc)=Z\rtimes\calG(\bar\cc)$ with $Z\isom\dbZ_p$.
The principal tuples of $\calG(\cc)$ are $\calA=(Z_1\nek Z_r,Z)$, where $(Z_1\nek Z_r$) is a principal tuple of $\calG(\bar\cc)$,  with the same root.

Note that always $Z_1\nek Z_r\leq\Ker(\theta)$, so indeed $Z_1\times\cdots\times Z_r\leq G$.

\begin{exam}
\label{example-of-compatibility-ex1}
\rm
Let $\cc=\langle\langle \calB_1\rangle\diamond\langle\calB_2\rangle\rangle$ for $\calB_1,\calB_2\in\grB$.
For convenience we write
\[
\calG(\cc)=Z\rtimes((Z_1\rtimes\calB_1)*(Z_2\rtimes\calB_2)),
\]
where $Z,Z_1,Z_2$ are copies of $\dbZ_p$.
There are two principal tuples $\calA_1,\calA_2$ in $\calG(\cc)$, namely, $(Z_i,Z)$ with the root $\calB_i$, $i=1,2$.
\end{exam}

More generally, let $\dd\leq\cc$ be elementary type constructions over $\grB$.
We will say that the principal tuple $\calA=(Z_1\nek Z_r)$ of $\calG(\cc)$ is \emph{compatible} with $\dd$ if, roughly speaking, its root already appears in $\dd$;
then we define a ``restricted" principal tuple $\calA_\dd=(Z_{i_1}\nek Z_{i_k})$ of $\calG(\dd)$, where $1\leq i_1<\cdots<i_k\leq r$, with the same root as $\calA$.
Specifically, the definition is by induction on the structure of $\dd$, as follows:

\medskip

$\bullet$ \quad
If $\dd=\calB$, where $\calB\in\grB$ is the root of $\calA$, then $\calA$ and $\dd$ are compatible, and we set $\calA_\dd$ to be the trivial $\dbZ_p$-tuple with the root $\calB$.

$\bullet$ \quad
If $\dd=\calB'$, where $\calB'\in\grB$ is not the root of $\calA$ (this includes the case where $\calB'$ and the root of $\calA$ are equal as pro-$p$ pairs, but have different  location in $\cc$), then $\calA$ and $\dd$ are not compatible.

$\bullet$ \quad
If $\dd=\dd_1\diamond\dd_2$ for nontrivial subconstructions $\dd_1,\dd_2$ of $\cc$, and $\calA$ is compatible with some $\dd_i$, $i=1,2$, then $\calA$ is compatible with $\dd$,  and $\calA_\dd$ is the image of $\calA_{\dd_i}$ under the embedding $\iota_{\dd_i,\dd}\colon\calG(\dd_i)\hookrightarrow\calG(\dd)$ with the same root.
Note that $i$ is determined by the location of the root, so $\calA$ cannot be compatible with both $\dd_1$ and $\dd_2$.

$\bullet$ \quad
If $\dd=\dd_1\diamond\dd_2$ for nontrivial subconstructions $\dd_1,\dd_2$ of $\cc$, and $\calA$ is compatible with neither $\dd_1$ nor $\dd_2$, then $\calA$ is not compatible with $\dd$.

$\bullet$ \quad
If $\dd=\langle\bar\dd\rangle$ for a subconstruction $\bar\dd$ of $\cc$, and $\calA$ is compatible with $\bar\dd$, then it is compatible with $\dd$.
We write $\calA_{\bar\dd}=(Z_{i_1}\nek Z_{i_{k-1}})$ and take $Z\isom\dbZ_p$ with $\calG(\dd)=Z\rtimes\calG(\bar\dd)$.
We set $\calA_\dd=(Z_{i_1}\nek Z_{i_{k-1}},Z)$ with the same root as $\calA_{\bar\dd}$, and where $Z_{i_1}\nek Z_{i_{k-1}}$ are identified with their images under $\iota_{\bar\dd,\dd}$.

$\bullet$ \quad
If $\dd=\langle\bar\dd\rangle$ for a subconstruction $\bar\dd$ of $\cc$, and $\calA$ is not compatible with $\bar\dd$, then $\calA$ is not compatible with $\dd$.

\begin{exam}
\label{example-of-compatibility-ex2}
\rm
Let $\cc$ and $\calA_1,\calA_2$ be as in Example \ref{example-of-compatibility-ex1}.
For $i=1,2$ consider the subconstruction $\dd_i=\langle\calB_i\rangle$ of $\cc$ with the associated pro-$p$ pair $\calG(\dd_i)=Z_i\rtimes\calB_i$, where $Z_i\isom\dbZ_p$.
Then, e.g., the principal tuple $\calA_1$ of $\calG(\cc)$ is compatible with $\dd_1$, but not with $\dd_2$ (even when $\calB_1=\calB_2$).
Further, $(\calA_1)_{\dd_1}=(Z_1)$ with the root $\calB_1$.
\end{exam}

\begin{rems}
\label{observations-about-compatibility-rems}
\rm
\begin{enumerate}
\item[(1)]
When $\dd=\cc$, every principal tuple $\calA$  for $\cc$ is compatible with $\cc$, and one has $\calA_\cc=\calA$.
\item[(2)]
If $\dd,\dd'$ are subconstructions of $\cc$ with $\dd\geq\dd'$, and if $\calA$ is compatible with $\dd'$, then it is also compatible with $\dd$.
\end{enumerate}
\end{rems}

The \emph{extension rank} $e(\cc)$ of an elementary type construction $\cc$ over $\grB$ is defined by
\[
e(\cc)=\sup_{\calA}\rk(\calA),
\]
where $\calA$ ranges of all principal tuples for $\cc$.

\begin{rems}
\label{properties-of-the-extension-rank-rems}
\rm
(1)\
The following rules follow from the definition of principal tuples:
\begin{enumerate}
\item[(i)]
If $\cc=\calB$, where $\calB\in\grB$, then $e(\cc)=0$;
\item[(ii)]
If $\cc=\cc_1\diamond\cc_2$ for nontrivial elementary type constructions $\cc_1,\cc_2$ over $\grB$, then $e(\cc)=\max(e(\cc_1),e(\cc_2))$;
\item[(iii)]
If $\cc=\langle\bar\cc\rangle$ for an elementary type construction $\bar\cc$ over $\grB$, then $e(\cc)=e(\bar\cc)+1$.
\end{enumerate}

(2) \quad
The extension rank $e(\cc)$ is not an invariant of $\calG(\cc)$.
For instance, if $p=2$, $\cc=\calE\diamond\calE$ and $\cc'=\langle\calE\rangle$, then the pro-$2$ pairs $\calG(\cc),\calG(\cc')$ are isomorphic (see Example \ref{c to G is not injective}), whereas $e(\cc)=0$ and $e(\cc')=1$.
\end{rems}

\section{Lifting of Automorphisms}
\label{section on automorphisms}

In the situation of Proposition \ref{normalization lemma}, when the $\dbZ_p$-tuple $\calA$ is a principal tuple of an elementary type pro-$p$ pair $\calG(\cc)$, we will need to know that the $\calA$-automorphism $\alp$ is induced by an automorphism of the pro-$p$ pair $\calG(\cc)$.
This will be a consequence of the following general result:

\begin{thm}
\label{automorphism theorem}
Let $\cc$ be an elementary type construction over $\grB$ and let $\calA=(Z_1\nek Z_r)$ be a principal tuple in $\calG(\cc)$.
Then every $\calA$-automorphism $\alpha$ of $Z_1\times\cdots\times Z_r$ extends to an automorphism of the pro-$p$ pair $\calG(\cc)$.
\end{thm}

The proof will be based on the following technical proposition.

\begin{prop}
\label{the morphism eta}
Let $\bar\cc,\bar\dd$ be elementary type constructions over $\grB$ such that $\bar\cc\geq\bar\dd$, and let $\cc=\langle\bar\cc\rangle$.
Let $\calA=(Z_1\nek Z_r)$ be a principal tuple in $\calG(\cc)$ of rank $r\geq1$, which is compatible with $\bar\dd$, and let
$\calA_{\bar\dd}=(Z_{i_1}\nek Z_{i_k})$ with $k<r$.
Let $\alpha$ be an $\calA$-automorphism such that $\alpha(z)=z\pmod{Z_r}$ for every $z\in Z_{i_j}$, $j=1,2\nek k$.
Then there is a morphism $\eta=\eta_{\bar\dd,\cc}\colon \calG(\bar\dd)\to\calG(\cc)$ which coincides with $\alp$ on $Z_{i_1}\nek Z_{i_k}$, and such that $\pi_{\cc,\bar\cc}\circ\eta=\iota_{\bar\dd,\bar\cc}$ (with $\pi_{\cc,\bar\cc}$ and $\iota_{\bar\dd,\bar\cc}$ as in Remark \ref{embedding of subconstruction}).
\end{prop}
\[
\xymatrix{
\calG(\cc)\ar[r]^{\pi_{\cc,\bar\cc}}&\calG(\bar\cc)\\
&\calG(\bar\dd)\ar[ul]^{\eta=\eta_{\bar\dd,\cc}}\ar[u]_{\iota_{\bar\dd,\bar\cc}}
}
\]
\begin{proof}
We abbreviate  $Z:=Z_r=\Ker(\pi_{\cc,\bar\cc})$, and write $\calG(\bar\cc)=(\bar G,\bar\theta)$ and $\calG(\cc)=Z\rtimes\calG(\bar\cc)=(G,\theta)$.
We argue by induction on the structure of $\bar\dd$ (for arbitrary $\cc$, $\bar\cc$, $\calA$ satisfying the above conditions).

\medskip

$\bullet$ \quad
If $\bar\dd=\calB\in\grB$, then $\calA_{\bar\dd}$ is the trivial $\dbZ_p$-tuple, so we may take $\eta=\iota_{\calB,\cc}$.

\medskip

$\bullet$ \quad
Suppose that $\bar\dd=\bar\dd_1\diamond\bar\dd_2$ for nontrivial elementary type constructions $\bar\dd_1,\bar\dd_2$.
We may assume that $\bar\calA$ is compatible with $\bar\dd_1$ but not with $\bar\dd_2$.
The induction hypothesis yields a morphism $\eta_1\colon\calG(\bar\dd_1)\to\calG(\cc)$ satisfying the required properties with $\bar\dd$ replaced by $\bar\dd_1$.
We take
\[
\eta=\eta_1*\iota_{\bar\dd_2,\cc}\colon \calG(\bar\dd)=\calG(\bar\dd_1)*\calG(\bar\dd_2)\to\calG(\cc).
\]

On the subgroups $Z_{i_1}\nek Z_{i_k}$ in $\calA_{\bar\dd}=\calA_{\bar\dd_1}$ we have $\eta=\eta_1=\alpha$.

Moreover,  on $\calG(\bar\dd_1)$ we have
\[
\pi_{\cc,\bar\cc}\circ\eta=\pi_{\cc,\bar\cc}\circ\eta_1=\iota_{\bar\dd_1,\bar\cc}=\iota_{\bar\dd,\bar\cc},
\]
and on $\calG(\bar\dd_2)$ we have
\[
\pi_{\cc,\bar\cc}\circ\eta=\pi_{\cc,\bar\cc}\circ\iota_{\bar\dd_2,\cc}=\iota_{\bar\dd_2,\bar\cc}=\iota_{\bar\dd,\bar\cc}.
\]
Therefore $\pi_{\cc,\bar\cc}\circ\eta=\iota_{\bar\dd,\bar\cc}$ on $\calG(\bar\dd)=\calG(\bar\dd_1)*\calG(\bar\dd_2)$.

\medskip

$\bullet$ \quad
Suppose that $\bar\dd=\langle\bar\dd'\rangle$ for an elementary type construction $\bar\dd'$ over $\grB$.
We write $\calG(\bar\dd')=(\bar G',\bar\theta')$.
Then  $\calG(\bar\dd)=Z'\rtimes\calG(\bar\dd')$ and $\calG(\langle\bar\dd\rangle)=Z''\rtimes\calG(\bar\dd)$, where $Z'=Z_{i_k}$ and $Z''\isom\dbZ_p$.
Further, $\calA$ is compatible with $\bar\dd'$, and we have
\[
\calA_{\langle\bar\dd\rangle}=(Z_{i_1}\nek Z_{i_{k-1}},Z',Z''), \quad
\calA_{\bar\dd'}=(Z_{i_1}\nek Z_{i_{k-1}}),
\]
 with the same root as $\calA$.
For convenience, we identify $Z''=Z$, so
\[
\calG(\langle\bar\dd\rangle)=Z\rtimes\calG(\bar\dd)=Z\rtimes(Z'\rtimes\calG(\bar\dd'))=\langle Z,Z'\rangle\rtimes\calG(\bar\dd').
\]

We apply the induction hypothesis with $\bar\dd$, $\bar\cc$, $\cc$, $\calA$ replaced by $\bar\dd'$, $\bar\dd=\langle\bar\dd'\rangle$, $\langle\bar\dd\rangle=\langle\langle\bar\dd'\rangle\rangle$, $\calA_{\langle\bar\dd\rangle}$, respectively.
It yields a morphism $\bar\eta'\colon\calG(\bar\dd')\to Z\rtimes(Z'\rtimes\calG(\bar\dd'))$ which coincides with $\alpha$ on $Z_{i_1}\nek Z_{i_{k-1}}$ and such that $\pi_{\langle\overline\dd\rangle,\overline\dd}\circ\bar\eta'=\iota_{\overline\dd',\overline\dd}$.
Thus there is a commutative diagram
\begin{equation}
\label{cd}
\xymatrix{
\calG(\cc)\ar[rr]^{\pi_{\cc,\bar\cc}}&&\calG(\bar\cc)\\
\calG(\langle\bar\dd\rangle)\ar[u]^{\iota_{\langle\overline\dd\rangle,\cc}}\ar[rr]^{\pi_{\langle\overline\dd\rangle,\overline\dd}}&&\calG(\bar\dd)\ar[u]_{\iota_{\overline\dd,\overline\cc}}\\
&&\calG(\bar\dd')\ar[llu]^{\bar\eta'}\ar[u]_{\iota_{\overline\dd',\overline\dd}}.
}
\end{equation}

The assumption on $\alp$ implies that it maps $Z'$ into $Z\times Z'=\langle Z,Z'\rangle$.
By Remark \ref{extended semidirect product relation} (applied with respect to the extension $\calG(\langle\bar\dd\rangle)=\langle Z,Z'\rangle\rtimes\calG(\bar\dd')$), for every $z'\in Z'$ and every $\bar g'\in\bar G'$ one therefore has
\[
\bar\eta'(\bar g')\cdot\alp(z')\cdot\bar\eta'(\bar g')\inv=\alp(z')^{\theta(\bar\eta'(\bar g'))}=\alp(z')^{\bar\theta'(\bar g')}.
\]
The universal property of extensions (see \S\ref{section on cyclotomic pro-p pairs}) yields a morphism
\[
\alp|_{Z'}\rtimes\bar\eta'\colon\calG(\bar\dd)=Z'\rtimes\calG(\bar\dd')\to\calG(\langle\bar\dd\rangle)
,
\]
which is $\alp$ on $Z'$ and $\bar\eta'$ on $\calG(\bar\dd')$.
We take the morphism $\eta\colon \calG(\bar\dd)\to\calG(\cc)$ to be the composition of this morphism with $\iota_{\langle\bar\dd'\rangle,\cc}$, namely,
\[
\calG(\bar\dd)=Z'\rtimes\calG(\bar\dd')\xrightarrow{\alp|_{Z'}\rtimes\bar\eta'}\calG(\langle\bar\dd\rangle)\xrightarrow{\iota_{\langle\bar\dd\rangle,\cc}}\calG(\cc).
\]
Note that it coincides with $\alp$ on $Z_1\nek Z_{i_{k-1}}\subseteq \bar G'$ as well as on $Z'=Z_{i_k}$.

By the commutativity of (\ref{cd}), on $\calG(\bar\dd')$ we have $\pi_{\cc,\bar\cc}\circ\eta=\iota_{\overline\dd,\overline\cc}\circ\iota_{\overline\dd',\overline\dd}=\iota_{\bar\dd',\bar\cc}=\iota_{\bar\dd,\bar\cc}$.

By the assumption on $\alpha$, we have $\pi_{\cc,\bar\cc}\circ\alp|_{Z'}=\id_{Z'}$, so $\pi_{\cc,\bar\cc}\circ\eta=\iota_{\bar\dd,\bar\cc}$ also on $Z'$, and therefore on all of $\calG(\bar\dd)=Z'\rtimes\calG(\bar\dd')$, as desired.
\end{proof}

\medskip

\begin{proof}[Proof of Theorem \ref{automorphism theorem}]
We proceed by induction on the structure of $\cc$.

\medskip

$\bullet$ \quad
If $\cc=\calB$ is a building block in $\grB$, then $\calA$ is the trivial tuple with root $\calB$, and we take $\gamma=\id_{\calB}$.

\smallskip

$\bullet$ \quad
Suppose that $\cc=\cc_1\diamond\cc_2$ for nontrivial elementary type constructions $\cc_1,\cc_2$ over $\grB$.
Without loss of generality, $\calA$ is a principal tuple in $\calG(\cc_1)$.
The inductive assumption yields an automorphism $\gamma_1$ of $\calG(\cc_1)$ which extends $\alp$.
Then the automorphism $\gamma=\gamma_1*\id_{\calG(\cc_2)}$ of $\calG(\cc)=\calG(\cc_1)*\calG(\cc_2)$ also extends $\alp$, as required.

\smallskip

$\bullet$ \quad
Suppose that $\cc=\langle\bar\cc\rangle$ for an elementary type construction $\bar\cc$ over $\grB$.
Thus $\calG(\cc)=Z_r\rtimes\calG(\bar\cc)$.
The principal tuple $\calA=(Z_1\nek Z_{r-1},Z_r)$ is compatible with $\bar\cc$, and $\calA_{\bar\cc}=(Z_1\nek Z_{r-1})$ is a principal tuple in $\calG(\bar\cc)$ with the same root as $\calA$.

Let $\bar\alp_{r-1}$ be, as in \S\ref{section on principal tuples}, the composition of the restriction $\alp|_{Z_1\times\cdots\times Z_{r-1}}$ with the projection  $Z_1\times\cdots\times Z_{r-1}\times Z_r\to Z_1\times\cdots\times Z_{r-1}$.
By Lemma \ref{induced automorphisms}, it is an $\calA_{\bar\cc}$-automorphism.
The inductive assumption on  $\bar\cc$ therefore yields an automorphism $\bar\gamma$ of the pro-$p$ pair $\calG(\bar\cc)$ which coincides with $\bar\alp_{r-1}$ on $Z_1\nek Z_{r-1}$.
Thus $\bar\gam(z)=\alp(z)\pmod{Z_r}$ for every $z\in Z_j$, $j=1,2\nek r-1$.

Furthermore, the restriction $\alpha|_{Z_r}$ is an automorphism of $Z_r$, by Lemma \ref{induced automorphisms} again.
Hence we have an automorphism of pro-$p$ pairs
\[
\gamma'=\alpha|_{Z_r}\rtimes\bar\gamma\colon \calG(\cc)=Z_r\rtimes\calG(\bar\cc)\to\calG(\cc)=Z_r\rtimes\calG(\bar\cc).
\]

Let $\alpha'$ be the restriction of $\gamma'$ to $Z_1\times\cdots\times Z_r$.
Then
\[
\alpha'|_{Z_r}=\alpha|_{Z_r}, \quad \alpha'|_{Z_j}=\bar\gamma|_{Z_j}=\bar\alp_{r-1}|_{Z_j}, \quad j=1,2\nek r-1.
\]
In particular, $\alpha'$ is a pro-$p$ group automorphism of $Z_1\times\cdots\times Z_{r-1}\times Z_r$.
By Lemma \ref{induced automorphisms}, $\alpha'$ is an $\calA$-automorphism.
Therefore  $\beta=(\alpha')^{-1}\circ\alpha$ is also an $\calA$-automorphism.
We have $\beta|_{Z_r}=\id_{Z_r}$, and $\beta(z)=z\pmod{Z_r}$ for $z\in Z_j$, $j=1,2\nek r-1$.

We now apply Proposition \ref{the morphism eta} with respect to the subconstruction $\bar\dd=\bar\cc$ and $\beta$.
We obtain a morphism of pro-$p$ pairs $\bar\eta\colon\calG(\bar\cc)\to\calG(\cc)$ which  coincides with $\beta$ on $Z_1\nek Z_{r-1}$, and such that $\pi_{\cc,\bar\cc}\circ\bar\eta=\id_{\calG(\bar\cc)}$.
By Lemma \ref{morphisms into extensions}, $\bar\eta$ extends to an automorphism $\eta$ of $\calG(\cc)$
such that $\eta|_{Z_r}=\id_{Z_r}$.
Then $\eta$ and $\beta$ coincide on $Z_1\nek Z_{r-1},Z_r$.
Setting $\gamma=\gamma'\circ\eta$, we obtain that
\[
\gamma|_{Z_1\times\cdots\times Z_r}=\gamma'|_{Z_1\times\cdots\times Z_r}\circ\beta=\alpha'\circ\beta=\alpha.
\qedhere
\]
\end{proof}

\bigskip

\section{Factoring of Epimorphisms}
\label{section on factoring of epimorphisms}

Let $\grB$ be again a nonempty class of building blocks.
Given an elementary type construction $\cc$ over $\grB$, we write $G(\cc)$ for the underlying pro-$p$ group of $\calG(\cc)$.

\begin{prop}
\label{quotient-by-standard-element-prop}
Let $\cc$ be an elementary type construction over $\grB$, and let $\calA=(Z_1\nek Z_r)$ be a principal tuple in $\calG(\cc)$.
Consider a pro-$p$ group $\bar G$ and a pro-$p$ group homomorphism $\rho\colon G(\cc)\to \bar G$.
Let $1\leq l\leq r$.
There exist a subconstruction $\cc'$ of $\cc$ and a pro-$p$ homomorphism $\rho'\colon G(\cc')\to \bar G$ with the following properties:
\begin{enumerate}
\item[(i)]
$\rho'\circ\pi_{\cc,\cc'}=\rho$ on  $G(\cc)$; and
\item[(ii)]
One has $\cc'\neq\cc$ if and only if $\rho(Z_l)=\{1\}$.
\end{enumerate}
\end{prop}

\begin{proof}
We argue by induction on the structure of $\cc$.

\smallskip

$\bullet$ \quad
If $\cc=\calB$ is a building block in $\grB$, then $\calG(\cc)=\calB$ and $r=0$, so there is nothing to prove.

\smallskip

$\bullet$ \quad
Suppose that $\cc=\cc_1\diamond\cc_2$ for nontrivial elementary type constructions $\cc_1,\cc_2$ over $\grB$.
Let $\rho_i$ be the restriction of $\rho$ to $G(\cc_i)$, considered as a closed subgroup of $G(\cc)$ via $\iota_{\cc_i,\cc}$, $i=1,2$.
Without loss of generality, $\calA$ is compatible with  $\cc_1$ and not with $\cc_2$.
Then $\calA_{\cc_1}=\calA$.
The induction hypothesis yields a subconstruction $\cc'_1$ of $\cc_1$  and a pro-$p$ group homomorphism $\rho'_1\colon G(\cc'_1)\to \bar G$,
such that (i) and (ii) hold with $\cc$, $\cc'$, $\rho$, $\rho'$ replaced by $\cc_1$, $\cc'_1$, $\rho_1$, $\rho'_1$, respectively.

If $\cc'_1\neq(1,1)$, then we take the subconstruction $\cc'=\cc'_1\diamond\cc_2$ of $\cc=\cc_1\diamond\cc_2$, and the pro-$p$ group homomorphism $\rho'=\rho'_1*_p\rho_2\colon G(\cc'_1)*_pG(\cc_2)\to \bar G$.
Note that here
\[
\pi_{\cc,\cc'}=\pi_{\cc_1,\cc'_1}*\id_{\calG(\cc_2)}\colon\calG(\cc)=\calG(\cc_1)*\calG(\cc_2)\to \calG(\cc')=\calG(\cc'_1)*\calG(\cc_2).
\]
Then (i) and (ii) hold.

If $\cc'_1=(1,1)$, then we take $\cc'=\cc_2$ and $\rho'=\rho_2$, and argue similarly.

\smallskip

$\bullet$ \quad
Suppose that $\cc=\langle\bar\cc\rangle$, where $\bar\cc$ is an elementary type construction over $\grB$ with $\calG(\bar\cc)=(G(\bar\cc),\bar\theta)$.
Then $\calG(\cc)=Z_r\rtimes\calG(\bar \cc)$ and $A_{\bar\cc}=(Z_1\nek Z_{r-1})$.
If $\rho(Z_l)\neq\{1\}$ then we just take $\cc'=\cc$ and $\rho'=\rho$.
We may therefore assume that $\rho(Z_l)=\{1\}$.

\case $l=r$. \rm
Then we take $\cc'=\bar\cc(\neq\cc)$, and note that $\rho$ factors as $\rho=\rho'\circ\pi_{\cc,\bar\cc}$ for some pro-$p$ group homomorphism $\rho'\colon G(\bar\cc)\to\bar G$.

\case $l\leq r-1$. \rm
Let $\bar\rho$ be the restriction of $\rho$ to $G(\bar\cc)$, considered as a closed subgroup of $G(\cc)$ via $\iota_{\bar\cc,\cc}$.
The induction hypothesis yields a subconstruction $\bar\cc'$ of $\bar\cc$ with $\bar\cc'\neq\bar\cc$
and a pro-$p$ group homomorphism $\bar\rho'\colon G(\bar\cc')\to \bar G$ such that $\bar\rho'\circ\pi_{\bar\cc,\bar\cc'}=\bar\rho$.
We put $\cc'=\langle\bar\cc'\rangle$ and view it as a proper subconstruction of $\cc=\langle\bar\cc\rangle$.
Note that
\[
\pi_{\cc,\cc'}=\id_{Z_r}\rtimes\pi_{\bar\cc,\bar\cc'}\colon \calG(\cc)=Z_r\rtimes\calG(\bar\cc)\to\calG(\cc')=Z_r\rtimes\calG(\bar\cc').
\]

\[
\xymatrix{
G(\cc)\ar[rrrr]^{\pi_{\cc,\bar\cc}}\ar@/_2pc/[rrdd]_{\rho}&&&&
   G(\cc')\ar@{-->}@/^2pc/[ddll]^{\rho'}\\
&G(\bar\cc)\ar[lu]_{\iota_{\bar\cc,\cc}}\ar[rr]^{\pi_{\bar\cc,\bar\cc'}}\ar[dr]^{\bar\rho}&& G(\bar\cc')\ar[ur]^{\iota_{\bar\cc',\cc'}}\ar[dl]_{\bar\rho'}\\
&&\bar G
}
\]
We define a pro-$p$ group homomorphism $\rho'\colon G(\cc')=Z_r\rtimes G(\bar\cc')\to \bar G$ by the rules $\rho'|_{Z_r}=\rho|_{Z_r}$ and $\rho'|_{G(\bar\cc')}=\bar\rho'$.
To see that $\rho'$ is well-defined, we need to verify that for $z\in Z_r$ and $\bar g'\in G(\bar\cc')$ one has
$\rho'(\bar g')\rho(z)\rho'(\bar g')\inv=\rho(z)^{\bar\theta'(\bar g')}$, where $\bar\theta'$ is the character of $\calG(\bar\cc')$.

To this end we choose $\bar g\in G(\bar\cc)$ with $\bar g'=\pi_{\bar\cc,\bar\cc'}(\bar g)$.
Then $\bar\rho'(\bar g')=\bar\rho(\bar g)$ and $\bar\theta(\bar g)=\bar\theta'(\bar g')$.
As $\bar gz\bar g\inv=z^{\bar\theta(\bar g)}$ in $\calG(\cc)$, we have
\[
\bar\rho'(\bar g')\rho(z)\bar\rho'(\bar g')\inv
=\bar\rho(\bar g)\rho(z)\bar\rho(\bar g)\inv=\rho(\bar gz\bar g\inv)
=\rho(z^{\bar\theta(\bar g)})=\rho(z)^{\bar\theta'(\bar g')},
\]
as desired.

Finally, on  $G(\bar\cc)$ one has $\rho'\circ\pi_{\cc,\cc'}=\bar\rho'\circ\pi_{\bar\cc,\bar\cc'}=\bar\rho=\rho$.
Furthermore, $\rho'\circ\pi_{\cc,\cc'}=\rho$ also on $Z_r$, and hence on $G(\cc)$.
\end{proof}

\begin{prop}
 \label{factorization-main-lemma}
Let $\cc$ be  an elementary type construction over $\grB$ with associated pro-$p$ pair $\calG(\cc)=(G,\theta)$,
let $\bar G$ be a pro-$p$ group, and let $\rho\colon G\to \bar G$ be a pro-$p$ group homomorphism.
Let  $\calA=(Z_1\nek Z_r)$ be a principal tuple in $\calG(\cc)$, and assume that  $\rho(V^k\calA)=\rho(V^{k+1}\calA)$ for some $1\leq k\leq r-1$ (with notation as in \S\ref{section on principal tuples}).

Then there exist a subconstruction $\cc'$ of $\cc$ with $\cc'\neq\cc$
and with associated pro-$p$ pair $\calG(\cc')=(G',\theta')$, as well as
an epimorphism $\phi\colon\calG(\cc)\to\calG(\cc')$, and a pro-$p$ group homomorphism $\rho'\colon G'\to \bar G$, such that $\rho=\rho'\circ\phi$ on $G$.
\end{prop}

\begin{proof}
We apply Proposition \ref{normalization lemma} for the restriction $\rho_0$ of $\rho$ to $A=Z_1\times\cdots\times Z_r$.
It yields an $\calA$-automorphism $\alpha$ of $A$ such that $\rho(\alpha(Z_k))=\{1\}$.
Theorem \ref{automorphism theorem} then yields an automorphism $\gamma$ of the pro-$p$ pair $\calG(\cc)$ which extends $\alpha$.
In particular, $\rho(\gamma(Z_k))=\{1\}$.
Proposition \ref{quotient-by-standard-element-prop} for the pro-$p$ group homomorphism $\rho\circ\gamma\colon G\to \bar G$
gives rise to a subconstruction $\cc'$ of $\cc$ with $\cc'\neq\cc$
and with associated pro-$p$ pair $\calG(\cc')=(G',\theta')$, and a pro-$p$ group homomorphism $\rho'\colon G'\to \bar G$,
such that $\rho\circ\gamma=\rho'\circ\pi_{\cc,\cc'}$ on $G$.
Then the epimorphism $\phi=\pi_{\cc,\cc'}\circ\gamma\inv\colon\calG(\cc)\to\calG(\cc')$ satisfies $\rho=\rho'\circ\phi$ on $G$, as desired.
\end{proof}

For a finite $p$-group $\bar G$, we set
\[
l(\bar G)=\max_E(\log_p|E|),
\]
where $E$ ranges over all abelian subgroups of $\bar G$.
Alternatively, $l(\bar G)$ is the maximal nonnegative integer $l$ such that there is a proper tower $E_0\supsetneq E_1\supsetneq\cdots\supsetneq E_l$ of abelian subgroups of $\bar G$.

\begin{prop} \label{factorization-proposition}
Let $\cc$ be  an elementary type construction over $\grB$ with associated pro-$p$ pair $\calG(\cc)=(G,\theta)$,
let $\bar G$ be a finite $p$-group, and $\rho\colon G\to \bar G$ a pro-$p$ group homomorphism.
Then there exist a subconstruction $\cc'$ of $\cc$ with extension rank $e(\cc')\le l(\bar G)$ and with associated pro-$p$ pair $\calG(\cc')=(G',\theta')$,
an epimorphism $\phi\colon\calG(\cc)\to\calG(\cc')$, and a pro-$p$ group homomorphism $\rho'\colon G'\to \bar G$ such that $\rho=\rho'\circ\phi$.
\end{prop}

\begin{proof}
We proceed by induction on the structure of $\cc$.
If $e(\cc)\le l(\bar G)$ (and in particular, if $\cc=\calB$ is in $\grB$), then we simply take $\cc'=\cc$, \ $\phi=\id_{\calG(\cc)}$, and $\rho'=\rho$.
We may therefore assume that $e(\cc)>l(\bar G)$.

By the definition of $e(\cc)$ (\S\ref{section on principal tuples}), there is a principal tuple $\calA=(Z_1\nek Z_r)$ in $\calG(\cc)$ of rank $r=e(\cc)$.
Let $A=Z_1\times\cdots\times Z_r$.
Since  $r>l(\bar G)\geq l(\rho(A))$,  the subgroups in the tower
\[
\rho(A)=\rho(V^0\calA)\supseteq\rho(V^1\calA)\supseteq\cdots\supseteq \rho(V^{r-1}\calA)\supseteq\rho(V^r\calA)=\{1\}
\]
cannot be all proper;
that is, there exists $0\le k\le r-1$ with $\rho(V^k\calA)=\rho(V^{k+1}\calA)$.

Therefore, Proposition \ref{factorization-main-lemma} yields a subconstruction $\cc''$ of $\cc$ with $\cc''\neq\cc$
and with associated pro-$p$ pair $\calG(\cc'')=(G'',\theta'')$,
as well as an epimorphism $\phi''\colon\calG(\cc)\to\calG(\cc'')$ and a pro-$p$ group homomorphism $\rho''\colon G''\to \bar G$ such that $\rho=\rho''\circ\phi''$.
The induction hypothesis, applied to $\cc''$ with this data, gives rise to a subconstruction $\cc'$ of $\cc''$ (whence also of $\cc$)
with associated pro-$p$ pair $\calG(\cc')=(G',\theta')$ and with $e(\cc')\le l(\bar G)$,
as well as to an epimorphism  $\phi'\colon\calG(\cc'')\to\calG(\cc')$, and a pro-$p$ group
homomorphism $\rho'\colon G'\to \bar G$ such that $\rho''=\rho'\circ\phi'$.
We take $\phi=\phi'\circ\phi''\colon\calG(\cc)\to\calG(\cc')$.
\end{proof}

\section{The Symbol Length under Elementary Operations}
\label{section on the symbol length under elementary operations}
Let $G$ be a pro-$p$ group and $n$ a positive integer.
Recall that a \emph{symbol} in $H^n(G,\dbF_p)$ is a cup product $\chi_1\cup\cdots\cup \chi_n$, where $\chi_1\nek \chi_n\in H^1(G,\dbF_p)$.
Given $\omega\in H^n(G,\dbF_p)$, its \emph{symbol length} $\syml(\omega)$ is the minimal integer $M\geq0$ such that $\omega$
can be presented as a sum of $M$ symbols in $H^n(G,\dbF_p)$.
If $\omega$ cannot be presented in this form, then we define $\syml(\omega)=\infty$.

In particular, one has $\syml(\omega)=0$ if and only if $\omega=0$.
When $n=1$ we have $\syml(\omega)=1$ for every $\omega\neq0$.

\begin{rem}
\label{basic properties of decomposition ranks}
\rm
\begin{enumerate}
\item[(1)]
For $\omega_1,\omega_2\in H^n(G,\dbF_p)$ one has
\[
 \syml(\omega_1+\omega_2)\le\syml(\omega_1)+\syml(\omega_2).
\]
\item[(2)]
For every $\omega_1\in H^{n_1}(G,\dbF_p)$ and $\omega_2\in H^{n_2}(G,\dbF_p)$ one has
\[
\syml(\omega_1\cup \omega_2)\leq\syml(\omega_1)\cdot\syml(\omega_2).
\]
\item[(3)]
For a pro-$p$ group homomorphism $\phi\colon G'\to G$ and $\omega\in H^n(G,\dbF_p)$ one has
\[
\syml(\phi^*(\omega))\le\syml(\omega),
\]
where $\phi^*\colon H^n(G,\dbF_p)\to H^n(G',\dbF_p)$ is the homomorphism induced by $\phi$.
\end{enumerate}
\end{rem}

Next we  track the behavior of the symbol length under the two operations on pro-$p$ pairs:
free pro-$p$ products and extensions.

First, let $G_1,G_2$ be pro-$p$ groups, considered as closed subgroups of $G=G_1*_pG_2$.
There are pullback (restriction) homomorphisms $\Res_{G_i}\colon H^n(G,\dbF_p)\to H^n(G_i,\dbF_p)$, $\omega\mapsto \Res_{G_i}(\omega)$, $i=1,2$.
The cohomology (graded) algebra $H^\bullet(G,\dbF_p)$ is the connected direct sum
\[
H^\bullet(G_1,\dbF_p)\sqcup H^\bullet(G_2,\dbF_p)
\]
of the cohomology algebras $H^\bullet(G_1,\dbF_p)$ and $H^\bullet(G_2,\dbF_p)$ via these homomorphisms
\cite{NeukirchSchmidtWingberg}*{Th.\ 4.1.4}.
We deduce:

\begin{lem}
\label{free-product-decomposition-ranks}
In this setup, let $\omega\in H^n(G,\dbF_p)$,  $n\geq1$.
Then
\[
\syml(\omega)=\max\bigl(\syml(\Res_{G_1}(\omega)),\syml(\Res_{G_2}(\omega))\bigr).
\]
\end{lem}

For a pro-$p$ group $K$, we denote
\[
M_n(K)=\sup\bigl\{\syml(\omega)\ \bigm|\ \omega\in H^n(K,\dbF_p)\bigr\}.
\]

\begin{cor}
\label{bounds on M for free products}
In the situation as above,
$M_n(G)=\max\bigl(M_n(G_1),M_n(G_2)\bigr)$.
\end{cor}

Next we consider a pro-$p$ pair $\bar\calG=(\bar G,\bar\theta)$ and the extension $\calG=(G,\theta)=\dbZ_p\rtimes\bar\calG$.
Thus $G=\dbZ_p\rtimes \bar G$.
Let $\pi=\pi_{\cc,\bar\cc}\colon G\to \bar G$ and $\iota=\iota_{\bar\cc,\cc}\colon\bar G\to G$ be as before.
Then $\Inf_G=\pi^*$ is the inflation map and $\Res_{\bar G}=\iota^*$ is the restriction map.
As $\pi\circ\iota=\id_{\bar G}$, we have $\Res_{\bar G}\circ\Inf_G=\id_{\bar G}$.

The structure of $H^\bullet(G,\dbF_p)$ was computed by Wadsworth \cite{Wadsworth83}*{\S3}, as follows:

For convenience we write $Z=\Ker(\pi)\cong\dbZ_p$.
Then $H^1(Z,\dbF_p)\cong\dbF_p$ and $H^m(Z,\dbF_p)=0$ for all $m\geq2$.
There exists $\beta\in H^1(G,\dbF_p)$ such that $\Res_Z(\beta)$ generates $H^1(Z,\dbF_p)$ and such that $\Res_{\bar G}(\beta)=0$ in $H^1(\bar G,\dbF_p)$.

We further claim that the standard action of $\bar G$ on $H^1(Z,\dbF_p)$ is trivial.
Recall that this action is given in general by $({}^{\bar g}\chi)(z)=g\chi(g^{-1}zg)$ for $\bar g\in\bar G$, $\chi\in H^1(Z,\dbF_p)$, and $z\in Z$, and where we take $g\in G$ with $\bar g=\pi(g)$.
In our case $\theta(g)\in 1+p\dbZ_p$, so $z^{\theta(g)}\cdot z^{-1}\in Z^p$.
Using Remark \ref{extended semidirect product relation} we obtain
\[
({}^{\bar g}\chi)(z)=g\chi(g^{-1}zg)=\chi(g^{-1}zg)=\chi(z^{\theta(g)})=\chi(z).
\]

Consequently, by \cite{Wadsworth83}*{Th.\ 3.1}, for every $n\geq2$ there is a short exact sequence
\begin{equation}
\label{exact sequence}
0\to H^n(\bar G,\dbF_p)\xrightarrow{\Inf_G}H^n(G,\dbF_p)\xrightarrow{t}H^{n-1}(\bar G,\dbF_p)\to0
\end{equation}
for a certain homomorphism $t$.
Moreover, this sequence is split by the section
\[
s\colon H^{n-1}(\bar G,\dbF_p)\to H^n(G,\dbF_p), \quad \bar\omega\mapsto\beta\cup\Inf(\bar\omega).
\]
As $\Res_{\bar G}(\beta)=0$, we have $\Res_{\bar G}\circ s=0$.

By (\ref{exact sequence}) and since $\Res_{\bar G}\circ\Inf_G=\id_{\bar G}$, we have $\Ker(t)\cap\Ker(\Res_{\bar G})=\{0\}$.
From $\Img(s)\subseteq\Ker(\Res_{\bar G})$, we therefore deduce that
$H^n(G,\dbF_p)=\Ker(t)\oplus\Ker(\Res_{\bar G})$ and $\Img(s)=\Ker(\Res_{\bar G})$.

\begin{lem}
\label{symbol length in extensions}
For every $\omega\in H^n(G,\dbF_p)$ with $n\geq2$, one has:
\[
\syml(\omega)\le \syml(\Res_{\bar G}(\omega))+M_{n-1}(\bar G).
\]
\end{lem}
\begin{proof}
In view of the previous comments,
\[
\omega-(\Inf_G\circ\Res_{\bar G})(\omega)\in\Ker(\Res_{\bar G})=\Img(s),
\]
i.e., there exists $\bar\omega\in H^{n-1}(G,\dbF_p)$ such that
\[
\omega=(\Inf_G\circ\Res_{\bar G})(\omega)+\beta\cup\Inf_G(\bar\omega).
\]
The assertion now follows from Remark \ref{basic properties of decomposition ranks}.
\end{proof}

\begin{cor}
\label{inequality for the bounds}
In the situation as above,
$M_n(G)\le M_n(\bar G)+M_{n-1}(\bar G)$.
\end{cor}

\section{Uniform Bounds}
\label{section on uniform bounds}

For our fixed nonempty class $\grB$ of building blocks, and for $m\geq1$, let
\[
\begin{split}
M_{m,0}^\grB&=\sup\{ M_m(B)\ |\ \calB=(B,\theta)\in\grB\} \\
&=\sup\bigl\{\syml(\omega) \bigm|\ \calB=(B,\theta)\in\grB, \ \omega\in H^m(B,\dbF_p)\bigr\}.
\end{split}
\]
For convenience, when $m=0$ we set $M_{0,0}^\grB=0$.
We define integers $M_m^\grB$ by $M_1^\grB=\max(M_{1,0}^\grB,1)$, and $M_m^\grB=M_{m,0}^\grB$ for $m\neq1$.
Thus $M_m^\grB=M_{m,0}^\grB$ unless $m=1$ and $\grB=\{(1,1)\}$.

Let $f(e,m)$ be the map uniquely defined for $e,m\geq0$ by the rules
\begin{equation}
\label{properties of f 1}
f(0,m)=M_m^\grB, \quad f(e,0)=0, \hbox{ and }
\end{equation}
\begin{equation}
\label{properties of f 2}
f(e,m)=f(e-1,m)+f(e-1,m-1) \hbox{ for } e,m\geq1.
\end{equation}
Explicitly,
\[
f(e,m)=\sum_{k=0}^{\min(e,m)}\binom ek M_{m-k}^\grB,
\]
where by convention, if $M_{m-k}^\grB=\infty$ for some $k$ in this sum, then $f(e,m)=\infty$.
Indeed, it is straightforward to verify that this map satisfies (\ref{properties of f 1}) and (\ref{properties of f 2}).
In particular, for $e\geq1$ we have
\[
f(e,1)=M_1^\grB, \quad
f(e,2)=M_2^\grB+eM_1^\grB.
\]
By (\ref{properties of f 2}), $f(e,m)\leq f(e',m)$ for $e\leq e'$.
Hence, for $e=\max(e_1,e_2)$ one has
\[
f(e,m)=\max(f(e_1,m),f(e_2,m)).
\]

\begin{lem}
\label{bound on MmG}
Let $m\geq1$ be an integer.
Let $\cc$ be an elementary type construction over $\grB$ with associated pro-$p$ pair $\calG(\cc)=(G,\theta)$ and extension rank $e=e(\cc)$.
Then  $M_m(G)\le f(e,m)$.
\end{lem}

\begin{proof}
We may assume that $f(e,m)<\infty$.
We argue again by induction on the construction of $\cc$.

\smallskip

$\bullet$ \quad
If $\cc=\calB$ is a building block in $\grB$, then $e=0$ and $M_m(G)\leq M_m^\grB=f(0,m)$.

\smallskip

$\bullet$ \quad
Suppose that $\cc=\cc_1\diamond\cc_2$ for nontrivial elementary type constructions $\cc_1,\cc_2$ over $\grB$, with $\calG(\cc_i)=(G_i,\theta_i)$ and $e_i=e(\cc_i)$, $i=1,2$.
By Remark \ref{properties-of-the-extension-rank-rems}(1), $e=e(\cc)=\max(e_1,e_2)$.
By induction, $M_m(G_i)\leq f(e_i,m)$, $i=1,2$.
In view of Corollary \ref{bounds on M for free products},
\[
M_m(G)=\max(M_m(G_1),M_m(G_2))\leq \max(f(e_1,m),f(e_2,m))=f(e,m).
\]

$\bullet$ \quad
If $\cc=\langle\bar\cc\rangle$ for an elementary type construction $\bar\cc$ over $\grB$, then $e\geq1$ and $e(\bar\cc)=e-1$ (Remark \ref{properties-of-the-extension-rank-rems}(1)).
We write $\calG(\bar\cc)=(\bar G,\bar\theta)$.
When $m\geq2$ we obtain using Corollary \ref{inequality for the bounds}, the induction hypothesis, and (\ref{properties of f 2}), that
\[
M_m(G)\le M_m(\bar G)+M_{m-1}(\bar G)  \le f(e-1,m)+f(e-1,m-1)=f(e,m).
\]
In the remaining case $m=1$ we have, since $G\neq\{1\}$ and by the definition of $M_1^\grB$,  that $M_1(G)=1\leq M_1^\grB=f(e,1)$.
\end{proof}

\medskip

Next we fix a pro-$p$ group $\bar G$, $n\geq2$, and a cohomology element $\bar\omega\in H^n(\bar G,\dbF_p)$.
From Lemma \ref{bound on MmG} we deduce:

\begin{cor}
\label{bound on siml of pullbacks}
Let $\cc$ be an elementary type construction over $\grB$ with associated pro-$p$ pair $\calG(\cc)=(G,\theta)$ and extension rank $e=e(\cc)$.
For every pro-$p$ group homomorphism $\rho\colon G\to \bar G$ one has
$\syml(\rho^*(\bar\omega))\le f(e,n)$.
\end{cor}

We can now prove Theorem A:

\begin{thm}
\label{uniform bound theorem}
Let $\bar G$ be a pro-$p$ group, let $n\geq2$, and let $\bar\omega\in H^n(\bar G,\dbF_p)$.
Assume that  $M_m^\grB<\infty$ for every $2\leq m\leq n$.
Then there is a positive integer $M=M(\bar G,n,\bar\omega)$ such that for every pro-$p$ pair $(G,\theta)$ of elementary type over $\grB$, and every pro-$p$ group homomorphism $\rho\colon G\to \bar G$, one has $\syml(\rho^*(\bar\omega))\le M$.

When $\bar G$ is finite, one can take $M=f(l(\bar G),n)$, where $l(\bar G)$ is as in \S\ref{section on factoring of epimorphisms}.
\end{thm}

\begin{proof}
We assume first that $\bar G$ is finite, so $l(\bar G)<\infty$.
The assumptions imply that $f(l(\bar G),n)<\infty$.

Let $\cc$ be an elementary type construction over $\grB$ such that $\calG(\cc)=(G,\theta)$.
Proposition \ref{factorization-proposition} yields a subconstruction $\cc'$ of $\cc$ with $e(\cc')\le l(\bar G)$ and with associated pro-$p$ pair $\calG(\cc')=(G',\theta')$, an epimorphism of pro-$p$ pairs $\phi\colon\calG(\cc)\to\calG(\cc')$, and a pro-$p$ group homomorphism
$\rho'\colon G'\to \bar G$ such that $\rho=\rho'\circ\phi$.
It follows from Remark \ref{basic properties of decomposition ranks}(3) and Corollary \ref{bound on siml of pullbacks} that
\[
\syml(\rho^*(\bar\omega))=\syml(\phi^*((\rho')^*(\bar\omega)))\le\syml((\rho')^*(\bar\omega))
\le f(e(\cc'),n)\le f(l(\bar G),n).
\]
Thus the assertion holds in this case with $M=f(l(\bar G),n)$.

In the general case, we have $H^*(\bar G,\dbF_p)=\varinjlim_{\bar N} H^*(\bar G/\bar N,\dbF_p)$, where $\bar N$ ranges
over the open normal subgroups of $\bar G$ \cite{NeukirchSchmidtWingberg}*{Prop.\ 1.2.5}.
Let $\phi_{\bar N}\colon \bar G\to\bar G/\bar N$ be the projection epimorphism.
Thus $\bar\omega=\phi_{\bar N}^*(\bar\omega_{\bar N})$ for some $\bar N$ and some $\bar\omega_{\bar N}\in H^n(\bar G/\bar N,\dbF_p)$.
We apply the finite case, with $\bar\omega$, $\rho$ replaced by $\bar\omega_{\bar N}$, $\phi_{\bar N}\circ\rho$, respectively, to obtain that
\[
\syml(\rho^*(\bar\omega))=\syml(\rho^*(\phi_{\bar N}^*(\bar\omega_{\bar N})))=\syml((\phi_{\bar N}\circ\rho)^*(\bar\omega_{\bar N}))\leq f(l(\bar G/\bar N),n).
\qedhere
\]
\end{proof}

\begin{exam}
\label{constants for standard building blocks}
\rm
The finiteness assumption in the statement of Theorem \ref{uniform bound theorem} holds in particular when $\grB$ is the class of standard building blocks (Example \ref{standard elementary type pairs example}).
Indeed, for $\calB=(B,\theta)\in\grB$ we compute:
\begin{enumerate}
\item[(1)]
If $\calB=(1,1)$, then $M_m(B)=0$ for every $m\geq2$.
\item[(2)]
If $\calB=(\dbZ_p,\theta)$, then $M_m(B)=0$ for $m\geq2$.
\item[(3)]
If $\calB=\calG_F(p)$ for a finite extension $F$ of $\dbQ_p$ which contains a root of unity of order $p$, then $B$ is a Demushkin group (see Remark \ref{Demushkin-groups-remark} or \cite{NeukirchSchmidtWingberg}*{Th.\ 7.5.11}).
Hence $M_2(B)=1$, and $M_m(B)=0$ for $m\geq3$.
\item[(4)]
If $p=2$ and $\calB=\calE$, then $B\cong\dbZ/2$.
We denote the nonzero element of $H^1(B,\dbF_2)$ by $\eps$.
For every $m\geq1$ one has $H^m(B,\dbF_2)=\{0,\eps\cup\cdots\cup\eps\}$, whence $M_m(B)=1$.
\end{enumerate}
Consequently, the finiteness assumption in the statement of Theorem \ref{uniform bound theorem} holds also for the classes of building blocks considered in Examples \ref{absolute elementary type pairs example}--\ref{pythagorean elementary type pairs example}.
Therefore Theorem  \ref{uniform bound theorem} applies to these cases as well.
\end{exam}

\section{Massey Products}
\label{section on Massey products}
In this final section we apply Theorem A to bound the symbol length of Massey product elements in $H^2(G,\dbF_p)$, where $G$ is the underlying group of a pro-$p$ pair $\calG=(G,\theta)$ of elementary type over the class $\grB$ of standard building blocks.
Assuming the Elementary Type Conjecture, this gives an upper bound on the symbol length of $m$-fold Massey product elements in $H^2(G_F(p),\dbF_p)$, for every field $F$ which contains a root of unity of order $p$ and such that $G_F(p)$ is a finitely generated pro-$p$ group (see Example \ref{standard elementary type pairs example}).
This bound is however not sharp.

Let $m\geq1$.
Recall that $\dbU_m(\dbF_p)$ is the group of all upper-triangular unipotent $(m+1)\times(m+1)$-matrices over $\dbF_p$.
It is a $p$-Sylow subgroup of $\GL_{m+1}(\dbF_p)$.

Our computation will be based on the following result due to Goozeff  \cite{Goozeff70}, for $p$ odd, and Barry \cite{Barry79}*{Th.\ 2.1} for arbitrary $p$:

\begin{thm}[Goozeff, Barry]
\label{Goozeff Barry}
The maximal order of an abelian $p$-subgroup of $\GL_{m+1}(\dbF_p)$, and hence of $\dbU_m(\dbF_p)$, is $p^{\lfloor(m+1)^2/4\rfloor}$.
Moreover, this order is attained by the elementary abelian $p$-group consisting of all matrices of the form
\[
\begin{bmatrix}
I_r&M\\
0&I_{m+1-r}
\end{bmatrix},
\]
where $r=\lfloor(m+1)/2\rfloor$, $M\in M_{r\times (m+1-r)}(\dbF_p)$, and $I_r,I_{m+1-r}$ denote the identity matrices.
\end{thm}
Note that indeed $r(m+1-r)=\lfloor(m+1)^2/4\rfloor$.

When $m\geq2$, let $\overline\dbU_m(\dbF_p)$ be, as before, the quotient of $\dbU_m(\dbF_p)$ by the subgroup consisting of all matrices $I+aE_{1,m+1}$ with $a\in \dbF_p$.
Alternatively, $\overline\dbU_m(\dbF_p)$ is the group of all \emph{partial} upper-triangular unipotent $(m+1)\times(m+1)$-matrices over $\dbF_p$ with the $(1,m+1)$-entry omitted.

\begin{lem}
\label{bound on l bar U}
One has $l(\overline\dbU_m(\dbF_p))\leq m^2/4+m-1$.
\end{lem}
\begin{proof}
There is an exact sequence
\[
0\to\dbF_p^{m-1}\to\overline\dbU_m(\dbF_p)\to\dbU_{m-1}(\dbF_p)\to1,
\]
where the right epimorphism is the projection on the upper-left $m\times m$-submatrix.
Let $E$ be an abelian subgroup of $\overline\dbU_m(\dbF_p)$, and let $\bar E$ be its image in $\dbU_{m-1}(\dbF_p)$.
Then $|E|\leq |\bar E|p^{m-1}$.
By Theorem \ref{Goozeff Barry}, $|\bar E|\leq p^{\lfloor m^2/4\rfloor}$, whence the assertion.
\end{proof}

We recall from the Introduction  the central extension of finite groups
\[
1\to\dbF_p\to\dbU_m(\dbF_p)\to\overline\dbU_m(\dbF_p)\to 1,
\]
which corresponds to a cohomology class $\bar\omega_m\in H^2(\overline\dbU_m(\dbF_p),\dbF_p)$ \cite{NeukirchSchmidtWingberg}*{Th.\ 1.2.4}.
We further recall that, for a pro-$p$ group $G$, the \emph{$m$-fold Massey product} $\langle\chi_1\nek\chi_m\rangle$ of cohomology elements $\chi_1\nek\chi_m\in H^1(G,\dbF_p)$ is the set of all pullbacks $\rho^*(\bar\omega_m)$ in $H^2(G,\dbF_p)$, where $\rho\colon G\to\overline\dbU_m(\dbF_p)$ is a pro-$p$ group homomorphism whose projections to the super-diagonal entries $(i,i+1)$ are $\chi_i$, $i=1,2\nek m$.

We now prove Theorem B:

\begin{thm}
\label{bound on Massey product elements}
Let $m\geq2$, let $(G,\theta)$ be a pro-$p$ pair of elementary type over the class $\grB$ of standard building blocks.
Let $\rho\colon G\to\overline\dbU_m(\dbF_p)$ be a pro-$p$ group homomorphism.
Then
\[
\syml(\rho^*(\bar\omega_m))\leq \dfrac {m^2}4+m.
\]
\end{thm}
\begin{proof}
By definition, $M_1^\grB\leq1$ and by Example \ref{constants for standard building blocks}, $M_2^\grB\leq1$.
Hence every $e\geq2$ one has $f(e,2)=M_2^\grB+eM_1^\grB\leq 1+e$.
By Theorem \ref{uniform bound theorem} and Lemma \ref{bound on l bar U},
\[
\syml(\rho^*(\bar\omega_m))\leq f(l(\overline\dbU_m(\dbF_p)),2)
\leq 1+l(\overline\dbU_m(\dbF_p))\leq
\dfrac{m^2}4+m.
\qedhere
\]
\end{proof}

\end{document}